\documentclass[11pt]{amsart}
\usepackage{geometry}                
\usepackage{graphicx}
\usepackage{amssymb}
\usepackage{epstopdf}
\usepackage{longtable}
\usepackage{psfrag}
\usepackage{overpic}
\usepackage{color}
\usepackage{caption2}
\usepackage{subfigure}

\title{A definite recursive relation and some statistical properties for M$\ddot{\mathrm{o}}$bius function}

\author{Rong Qiang Wei}
\address{College of Earth Sciences, University of Chinese Academy of Sciences, Beijing, PRC, 100049}
\email{wrq1973@ucas.ac.cn}
\date{}

\begin{document}
\maketitle
\begin{abstract}
  An elementary definite recursive relation for M$\ddot{\mathrm{o}}$bius function $\mu
  (n)$ is introduced by two simple ways. With this recursive relation, $\mu (n)$ can be
  calculated without directly knowing the factorization of the $n$. $\mu (1)
  \sim \mu (2 \times 10^7) $ are calculated recursively one by one. Based on
  these $2\times 10^7$ samples, the empirical probabilities of $\mu (n)$ of taking  
  $- 1,0,\mbox{\ and\ } 1$ in classic statistics are calculated and compared with the theoretical probabilities in
  number theory. The numerical consistency between these two kinds of
  probability show that $\mu (n)$ could be seen as an independent random
  sequence when $n$ is large. The expectation and variance of the $\mu
  (n)$ are $0$ and $6 n/ \pi^2$, respectively. Furthermore, we show that
  any conjecture of the Mertens type is false in probability sense, and
  present an upper bound for cumulative sums of $\mu (n)$ with a certain
  probability.
\end{abstract}

{\hspace{2.2em}\small Keywords:}

{\hspace{2.2em}\tiny M$\ddot{\mathrm{o}}$bius function  Recursive relation Probability sense Mertens conjecture}

\section{Introduction}

The M$\ddot{\mathrm{o}}$bius function is defined for a positive integer $n$ by

\begin{equation}\label{eq1}
\mu(n)  = \left\{ {\begin{array}{*{20}{c}}
1\\
0\\
{{{( - 1)}^k}}
\end{array}\begin{array}{*{20}{c}}
{n = 1{\mbox{\hspace{14em} }}}\\
{{\mbox{if \ }} n {\mbox{\ is divisible by a prime square\hspace{2em}}}}\\
{{\mbox{if\  }} n {\mbox{\  is the product of }}k{\mbox{ distinct primes}}}
\end{array}} \right.
\end{equation}

It is shown that M$\ddot{\mathrm{o}}$bius function and its associated
M$\ddot{\mathrm{o}}$bius transform are important for solving different
mathematical and/or scientific problems (eg., Schroeder, 2008). In physics, the M$\ddot{\mathrm{o}}$bius
 function and its associated M$\ddot{\mathrm{o}}$bius
transform are used in inverse black body radiation problem (eg., Chen, 1987; 1990),
inversion of specific heat data for phonon densities of states (eg., Chen et al., 1990), solution
of integral equations regarding Fermi and Bose systems (eg., Chen, 2010), inverse 
transmissivity problem (Ji et al., 2006), and so on. All of these studies 
are related to how to calculate $\mu(n)$ if special methods are not used.

To calculate the M$\ddot{\mathrm{o}}$bius function, many algorithms are
presented and most of them base on the factorization of its argument. A famous
one is vectorized sieving (eg., Lioen and Lune, 1994; Kuznetsov, 2011). On the other hand, 
in their book, Hardy and Wright (2008) showed that the M$\ddot{\mathrm{o}}$bius function is the sum of
the primitive $n$-th roots of unity, ie.,

\begin{equation}\label{eq2}
\mu (n) = \sum\limits_{\begin{array}{*{20}{c}}
{1 \le k \le n}\\
{\gcd (k,n) = 1}
\end{array}} {\exp (\frac{{2\pi ik}}{n})} 
\end{equation}

Formula (\ref{eq2}) can be used to calculate the M$\ddot{\mathrm{o}}$bius function
without knowing the factorization of $n$. However, the computational complexity is not low.

Here we introduce a definite recursive relation to calculate the M$\ddot{\mathrm{o}}$bius
function without directly knowing the factorization of $n$ as Formula
(\ref{eq2}) does, but the computational complexity is less. We calculate
M$\ddot{\mathrm{o}}$bius function from $\mu (1) $ to $\mu (2 \times 10^7)$ with
this recursive relation and discuss some statistical properties of the M$\ddot{\mathrm{o}}$bius function.

\section{A definite recursive relation for M$\ddot{\mathrm{o}}$bius function}\label{sec2}

A definite recursive relation for M$\ddot{\mathrm{o}}$bius function can be
introduced by two simple ways. One is from M$\ddot{\mathrm{o}}$bius transform, and the
other is from the Redheffer Matrix related to Mertens function which is the
cumulative sum of the M$\ddot{\mathrm{o}}$bius function.  The more general relation
for (poset) M$\ddot{\mathrm{o}}$bius function can be found in 
the Incidence Algebra (eg., https://en.wikipedia.org
/wiki/Incidence\_algebra).

\subsection{The recursive relation from M$\ddot{\mathrm{o}}$bius transfrom}

According to pair potential model for cohesive energy (Chen, 1994), the
cohesive energy $E$ for each atom in an infinite linear chain can be expressed
as a sum of pairwise potentials,

\begin{equation}\label{eq3}
E(x) = \sum\limits_{n = 1}^\infty  {\Phi (nx)}
\end{equation}

Using Chen-M$\ddot{\mathrm{o}}$bius formula (e.g., Chen, 2010; Wang, 2013),

\begin{equation}\label{eq4}
\Phi (x) = \sum\limits_{n = 1}^\infty  {\mu (n)E(nx)}
\end{equation}

We can write the following Matrix equality according to expression (\ref{eq3}) and
(\ref{eq4}),

\begin{equation}\label{eq5}
\left( {\begin{array}{*{20}{c}}
{E(x)}\\
{E(2x)}\\
{E(3x)}\\
{E(4x)}\\
{E(5x)}\\
{E(6x)}\\
 \vdots 
\end{array}} \right) = \left( {\begin{array}{*{20}{c}}
1&1&1&1&1&1& \cdots \\
0&1&0&1&0&1& \cdots \\
0&0&1&0&0&1& \cdots \\
0&0&0&1&0&0& \cdots \\
0&0&0&0&1&0& \cdots \\
0&0&0&0&0&1& \cdots \\
 \vdots & \vdots & \vdots & \vdots & \vdots & \vdots & \ddots 
\end{array}} \right)\left( {\begin{array}{*{20}{c}}
{\Phi (x)}\\
{\Phi (2x)}\\
{\Phi (3x)}\\
{\Phi (4x)}\\
{\Phi (5x)}\\
{\Phi (6x)}\\
 \vdots 
\end{array}} \right)
\end{equation}

\begin{equation}\label{eq6}
\left( {\begin{array}{*{20}{c}}
{\Phi (x)}\\
{\Phi (2x)}\\
{\Phi (3x)}\\
{\Phi (4x)}\\
{\Phi (5x)}\\
{\Phi (6x)}\\
 \vdots 
\end{array}} \right) = \left( {\begin{array}{*{20}{c}}
{\mu (1)}&{\mu (2)}&{\mu (3)}&{\mu (4)}&{\mu (5)}&{\mu (6)}& \cdots \\
0&{\mu (1)}&0&{\mu (2)}&0&{\mu (3)}& \cdots \\
0&0&{\mu (1)}&0&0&{\mu (2)}& \cdots \\
0&0&0&{\mu (1)}&0&0& \cdots \\
0&0&0&0&{\mu (1)}&0& \cdots \\
0&0&0&0&0&{\mu (1)}& \cdots \\
 \vdots & \vdots & \vdots & \vdots & \vdots & \vdots & \ddots 
\end{array}} \right)\left( {\begin{array}{*{20}{c}}
{E(x)}\\
{E(2x)}\\
{E(3x)}\\
{E(4x)}\\
{E(5x)}\\
{E(6x)}\\
 \vdots 
\end{array}} \right)
\end{equation}

Let $$\Phi = \left[ \Phi (x) \quad \Phi (2 x) \quad \Phi (3 x) \quad \Phi (4
x) \quad \Phi (5 x) \quad \Phi (6 x) \quad \ldots \right]^{\mbox{\tiny T}}$$ 
and $$E = \left[ E (x) \quad E (2 x) \quad E (3 x) \quad E (4 x)
\quad E (5 x) \quad E (6 x) \quad \ldots \right]^{\mbox{\tiny T}}$$

Matrix equality (\ref{eq5}) and (\ref{eq6}) can be rewritten as the following,

\begin{equation}\label{eq7}
E = U \Phi
\end{equation}

\begin{equation}\label{eq8}
\Phi = V E
\end{equation}

Obviously,

\begin{equation}\label{eq9}
U V = V U = V^T U^T = I
\end{equation}

Hence the values of the M$\ddot{\mathrm{o}}$bius function, which are the
elements of the first row of the matrix $V$, can be obtained from the inverse
matrix of $U$. Because the matrix $U$ is a triangular one, in which $U = \{u_{i j} \}$ with
 $u_{i j} = 1$ if and only if $i \vert j $, one can get,
 
 \begin{equation}\label{eq10}
 {v_{1i}} =  - \sum\limits_{k = 1}^{i - 1} {v_{1k}}{{u_{ki}}} ,i = 2,3, \cdots 
 \end{equation}

Based on recursive relation (\ref{eq10}), we can obtain the recursive relation for M$\ddot{\mathrm{o}}$bius
function as the following,

\begin{equation}\label{eq11}
\mu (n) =  - \sum\limits_{k = 1}^{n -1} {{l_{nk}}} \mu (k),n = 2,3, \cdots ;{l_{nk}} = \left\{ {\begin{array}{*{20}{c}}
1&{k|n}\\
0&{{\rm{else}}}
\end{array}} \right.
\end{equation}

\subsection{The recursive relation from Redheffer Matrix}

It is well known that the Mertens function, which is the cumulative sum of the
M$\ddot{\mathrm{o}}$bius function, is the determinant of the Redheffer matrix.
The Redheffer matrix $R = \{ r_{i j} \}$ is defined by $r_{i j} = 1$ if $j =1$ or 
$i \vert j $, and $r_{i j} = 0$ otherwise, ie.,

\begin{equation}\label{eq12}
R=\left( {\begin{array}{*{20}{c}}
1&1&1&1&1&1& \cdots \\
1&1&0&1&0&1& \cdots \\
1&0&1&0&0&1& \cdots \\
1&0&0&1&0&0& \cdots \\
1&0&0&0&1&0& \cdots \\
1&0&0&0&0&1& \cdots \\
 \vdots & \vdots & \vdots & \vdots & \vdots & \vdots & \ddots 
\end{array}} \right)
\end{equation}

$R$ can be decomposed as follows:

\begin{equation}\label{eq13}
R = \left( {\begin{array}{*{20}{c}}
0&0&0&0&0&0& \cdots \\
1&0&0&0&0&0& \cdots \\
1&0&0&0&0&0& \cdots \\
1&0&0&0&0&0& \cdots \\
1&0&0&0&0&0& \cdots \\
1&0&0&0&0&0& \cdots \\
 \vdots & \vdots & \vdots & \vdots & \vdots & \vdots & \ddots 
\end{array}} \right) + \left( {\begin{array}{*{20}{c}}
1&1&1&1&1&1& \cdots \\
0&1&0&1&0&1& \cdots \\
0&0&1&0&0&1& \cdots \\
0&0&0&1&0&0& \cdots \\
0&0&0&0&1&0& \cdots \\
0&0&0&0&0&1& \cdots \\
 \vdots & \vdots & \vdots & \vdots & \vdots & \vdots & \ddots 
\end{array}} \right) = S + U
\end{equation}

where $S = \{ s_{i j} \} = 1$ if and only if $j = 1$ and $i \neq 1$; $U$ is
the same to matrix equality (\ref{eq7}) with $u_{i j} = 1$ if and only if $i \vert j$.

It can be shown that the inverse of $U$ is $V$ which is in matrix equality
(\ref{eq8}), that is,

\begin{equation}\label{eq14}
V = \left\{ {{v_{ij}}} \right\} = \left\{ {\begin{array}{*{20}{c}}
{\mu (\frac{j}{i})}&{i|j}\\
0&{{\rm{else}}}
\end{array}} \right.
\end{equation}

In fact, the $i j $-th entry of the product of $U \times V$, $p_{i j}$, is,

\begin{equation}\label{eq15}
{p_{ij}} = \sum\limits_{k = 1}^n {{u_{ik}}} {v_{kj}}
\end{equation}

According to the definition of $U$ and $V$ in (14), $u_{i k} v_{k j}$ is 0
unless $i \vert k$ and $k \vert j $, which means that $p_{i j} =
0$ if $i \nmid j$ . If $i \vert j $, by using the well known
$$\sum\limits_{i|n} {\mu (i) = \left\{ {\begin{array}{*{20}{c}}
1&{n = 1}\\
0&{{\rm{else}}}
\end{array}} \right.}$$

one can get,
$${p_{ij}} = \sum\limits_{k(i|k{\rm{\  and \  }}k|j)} {\mu (\frac{j}{k}} ) = \sum\limits_{k'|(j/i)} {\mu (\frac{{j/i}}{{k'}}} ) = \left\{ {\begin{array}{*{20}{c}}
1&{j = i}\\
0&{{\rm{else}}}
\end{array}} \right.$$

Therefore, $U \times V = \{ p_{i j} \} = I$ and equality (14) holds. With the
same procedures in the subsection above, we can obtain the recursive relation
(\ref{eq11}) for the M$\ddot{\mathrm{o}}$bius function.

\section{Some statistical properties for M$\ddot{\mathrm{o}}$bius function}

With the recursive relation (\ref{eq11}), we calculated the
M$\ddot{\mathrm{o}}$bius function from $\mu (1)$ to $\mu (2 \times 10^7)$.
These values are used for the numerical test on some statistical properties of
M$\ddot{\mathrm{o}}$bius sequence $\mu (n)$, if $\mu (n)$ is seen as an
independent random sequence although it has a deterministic recursive rule. In
fact, as $n$ is large enough, the random assumption above is reasonable.

\subsection{The expectation and variance of the $\mu (n)$}\label{subsec3_1}

Firstly we calculated the probabilities of $\mu (n)$ of taking the values $-1$, 
0 and 1. In this respect, there are two useful results from Hardy and
Wright (2008) as follows,

1. $\mu (n) = \pm 1$ or $\vert \mu (n) \vert = 1$ if a number $n$ is squarefree, and
the probability ($p_t$) that a number should be squarefree is
$\frac{6}{\pi^2}$, more precisely,

\begin{equation}\label{eq16}
\overset{x}{\underset{n = 1}{\sum}} \vert \mu (n) \vert =
\frac{6}{\pi^2} x + O \left( \sqrt{x} \right)
\end{equation}

2. Among the squarefree numbers, those for $\mu (n) = 1$ and those for $\mu
(n) = - 1$ occur with about the same frequency.

Therefore, if $\mu (k)  (k = 1, 2, \ldots n)$ denotes the $k$-th value of $\mu(n)$ and $p_{t_{_k}}$ the probability, the corresponding distribution rule is shown in Table \ref{tb1}.

\begin{table}[htdp]
\caption{The distribution rule for $\mu(k)$}
\begin{center}
\begin{tabular}{c c c c}
\hline
$\mu (k)$ &$ -1$ &$0$ &$1$\\
\hline
$p_{t_{_k}}$&$3/\pi^2$&$1-6/\pi^2$&$3/\pi^2$\\
\hline
\end{tabular}
\end{center}
\label{tb1}
\end{table}%

The $p_t$ is different from that in the classic statistics (See
details in Hardy and Wright (2008), P. 354).  To use the methods in classic statistic,
it is necessary firstly to test the consistency between them. The classic
probability here is,

\begin{equation}\label{eq17}
p_e = \frac{{{N_{\mu (n) = m}}}}{n},{\rm{ }}(m =  - 1,0,1)
\end{equation}

where $p_e$ is the classic probability that $\mu (n) = m (m = - 1, 0, 1)$,
$N_{\mu (n) = m} $ is the frequency of $\mu (n) = m$. 

\begin{figure}[htb]
\setlength{\belowcaptionskip}{0pt}
\centering
\begin{overpic}[scale=0.6]{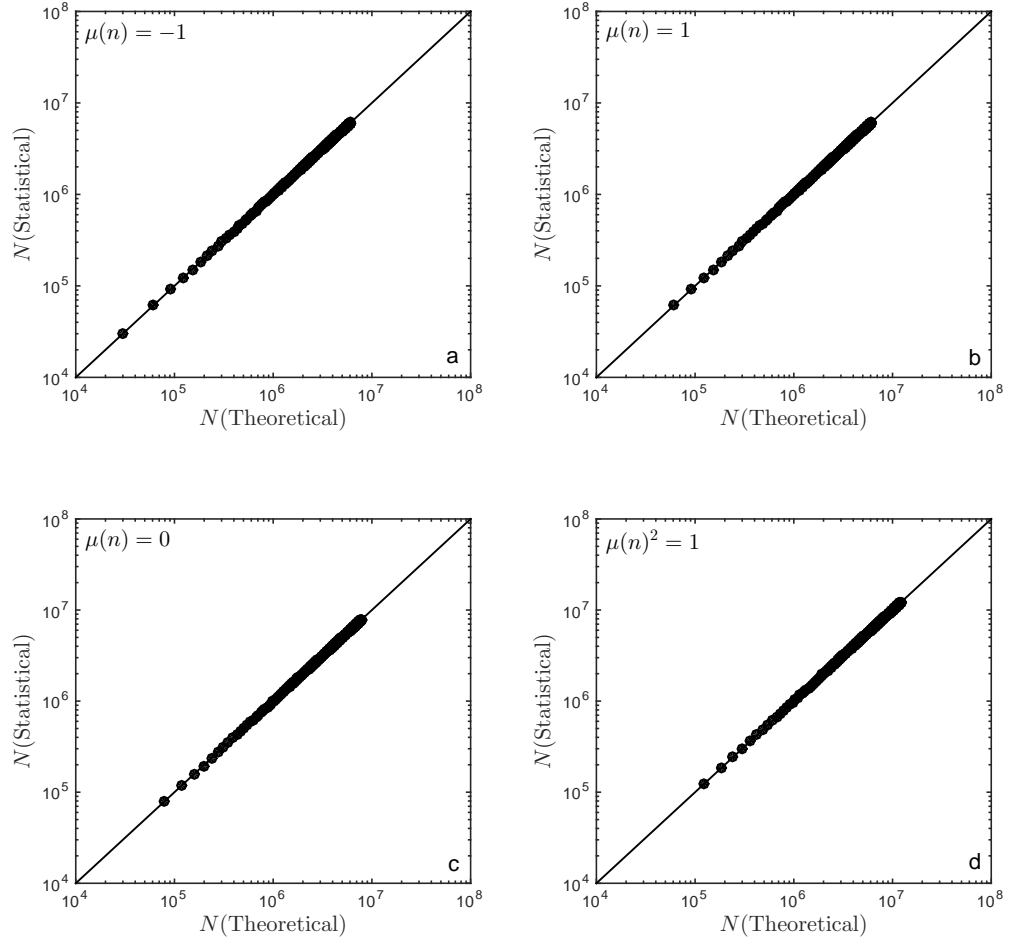}
\end{overpic}
\renewcommand{\figurename}{Fig.}
\caption{The comparison of frequencies observed ($N(\mbox{Statistical})$) with those calculated by number theory ($N(\mbox{Theorical})$) in blocks with different length of $N$. The solid line in each subfigure is the reference line. The related data are shown in Table \ref{tb3} and \ref{tb4} in the appendix.}
\label{fig1}
\end{figure}

\begin{figure}[htb]
\setlength{\belowcaptionskip}{0pt}
\centering
\begin{overpic}[scale=0.6]{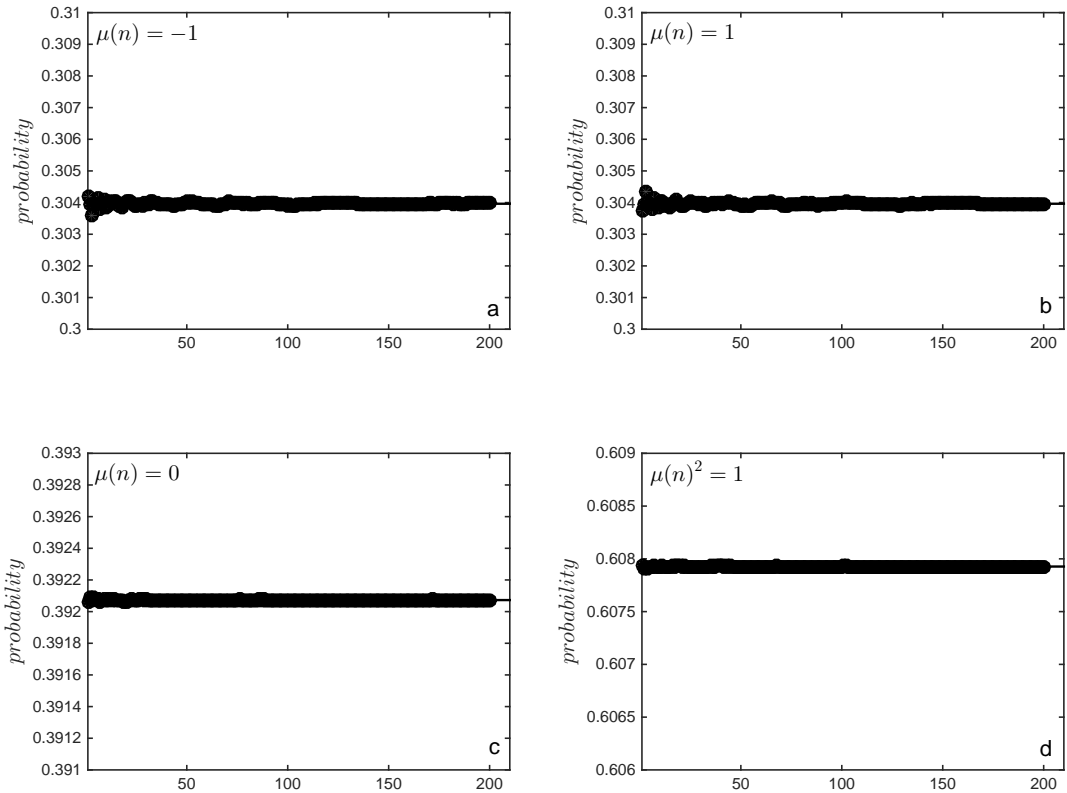}
\end{overpic}
\renewcommand{\figurename}{Fig.}
\caption{The comparison of classic probability $p_e$ (solid dotts) with the $p_t$ (solid lines) in blocks with different length of $N$. In each subfigure, the x-axis is ordinal number of each block. The related data are shown in Table \ref{tb3} and \ref{tb4} in the appendix.}
\label{fig2}
\end{figure}

We calculate the frequencies and $p_e$s of $\mu(n)$ of taking $-1,0,1$ and that of $\vert\mu(n)\vert=\mu^2(n)=1$ in 200 blocks with different length by using $2\times 10^7$ $\mu(n)$s above, respectively. 
Figure \ref{fig1} shows the comparison of these frequencies observed with those calculated by $N \times p_t$ in different blocks of length $N$. It can be seen that the frequencies observed are
consistent with those calculated. Figure \ref{fig2} shows the numerical comparison of
classic probability $p_e$  with the $p_t$.
It also can be observed that these two kinds of probabilities are numerically consistent.
Detail numerical results are shown in Table 3 and 4 in the appendix. These consistencies above show
that the $p_t$ is equivalent numerically to the classic
probability $p_e$ as defined in (\ref{eq17}). Similar numerical support can be found in
Good and Churchhouse (1968). Based on these numerical results, we can take $\mu(n)$ as an
independent random sequence although it has a deterministic recursive rule and use
classic statistical method to study it.

Accordingly, the expectation and variance of the $\mu (k)$ are $E (\mu (k)) = 0$
and $D (\mu (k)) = 6 / \pi^2$ from Table \ref{tb1}, respectively. These results are consistent with
the conjecture of Good and Churchhouse (1968). The expectation and variance of
the $\mu (k)$ will be used in the following section.

\subsection{Mertens conjecture in a statistical point of view}

Mertens function of a positive integer $n$ is defined as the cumulative
sums of $\mu (n)$,

\begin{equation}\label{eq18}
M (n) = \overset{n}{\underset{k = 1}{\sum}} \mu(k)
\end{equation}

An old conjecture, "Mertens conjecture" , proposed that $\vert M (n) \vert < n^{1/ 2}$ for all $n$. 
This was disproved by Odlyzko and te Riele (1985). In this subsection, we recheck Mertens conjecture in a statistical
point of view, for $\mu (n)$ is seen as an independent random sequence although it has a deterministic recursive rule.

According to central limit theorem, if $n$ is large enough, for any $x$, we have,

\begin{equation}\label{eq19}
\begin{array}{*{20}{c}}
{\mathop {\lim }\limits_{n \to \infty } P\left\{ {\frac{{\sum\limits_{k = 1}^n {\mu (k) - E(\sum\limits_{k = 1}^n {\mu (k)} )} }}{{\sqrt {D(\sum\limits_{k = 1}^n {\mu (k)} )} }} \le x} \right\}}\\
{}\\
{}
\end{array}\begin{array}{*{20}{l}}
{ = \mathop {\lim }\limits_{n \to \infty } P\left\{ {\frac{{M(n)}}{{\sqrt {6n/{\pi ^2}} }} \le x} \right\}}\\
{ = \int\limits_{ - \infty }^x {\frac{1}{{\sqrt {2\pi } }}\exp ( - \frac{{{t^2}}}{2}){\rm{d}}t}}\\
{ = \Phi (x)}
\end{array}
\end{equation}

where, $E \left( \sum^n_{k = 1} \mu (k) \right) = 0$, $D \left( \sum^n_{k = 1}
\mu (k) \right) = 6 n / \pi^2$ according to subsection \ref{subsec3_1}. And (\ref{eq19}) means,

\begin{equation}\label{eq20}
\frac{M (n)}{\sqrt{6 n / \pi^2}} \backsim N (0,1)
\end{equation}

as $n$ is large.

\begin{figure}[htb]
\setlength{\belowcaptionskip}{0pt}
\centering
\begin{overpic}[scale=0.6]{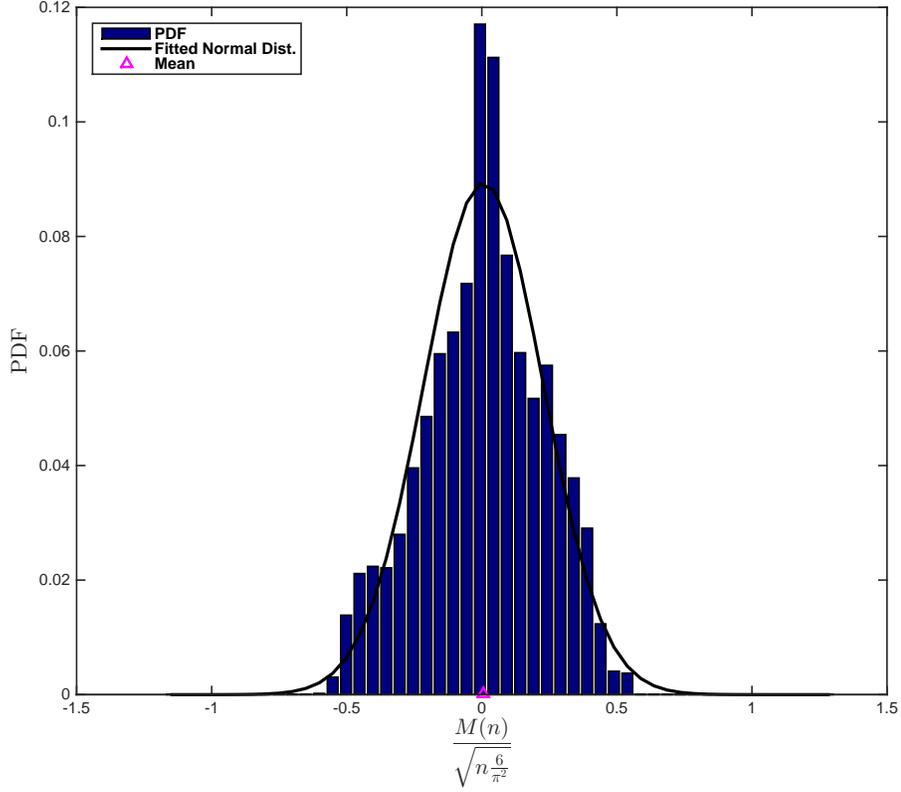}
\end{overpic}
\renewcommand{\figurename}{Fig.}
\caption{The probability density function for $\frac{M (n)}{\sqrt{6 n /\pi^2}}$ \ when $n = 500000$.}
\label{fig3}
\end{figure}

Figure \ref{fig3} shows the probability density function for $\frac{M (n)}{\sqrt{6 n /
\pi^2}}$ \ when $n = 500000$. It can be observed that the distribution of (\ref{eq20}) is
reasonable. Another similar numerical support for this can be found in Good and Churchhouse (1968).

With equality (\ref{eq19}), the probability of $M (n) > \sqrt{n}$ can be obtained.
Clearly,

\begin{equation}\label{eq21}
P\left\{ {M(n) > \sqrt n } \right\} = 1 - \Phi (\frac{1}{{\sqrt {6/{\pi ^2}} }}) \approx 0.0998
\end{equation}

That is, the probability of $M (n) > \sqrt{n}$ is about $0.0998$ but not $0$,
which means that Mertens conjecture is not true. Furthermore, any conjecture
of the Mertens type, viz.

\begin{equation}\label{eq22}
\left| {M(n)} \right| < C\sqrt n 
\end{equation}

where $C$ is any positive constant, is false, unless $C$ is large enough.

\subsection{Upper bound of cumulative sums of $\mu (n)$ sequence}

From (\ref{eq19}), one can get when $n$ is large,

\begin{equation}\label{eq23}
\frac{{\sum\limits_{k = 1}^n {\mu (k) - E(\sum\limits_{k = 1}^n {\mu (k)} )} }}{{\sqrt {D(\sum\limits_{k = 1}^n {\mu (k)} )} }} = \frac{{\sum\limits_{k = 1}^n {\mu (k) - nu} }}{{\sqrt n \sigma }} \sim N(0,1)
\end{equation}

Then a confidence interval for expectation $u$ with a known standard variance $\sigma$ and a probability of $1 - \alpha$ is,

\begin{equation}\label{eq24}
\left[ { - \frac{\sigma }{{\sqrt n }}{K_{\alpha /2}},{\rm{ }}\frac{\sigma }{{\sqrt n }}{K_{\alpha /2}}} \right]
\end{equation}

where,

\begin{equation}\label{eq25}
\int\limits_{ - {K_{\alpha /2}}}^{{K_{\alpha /2}}} {\frac{1}{{\sqrt {2\pi } }}} \exp ( - \frac{{{t^2}}}{2}){\rm{d}}t = 1 - \alpha 
\end{equation}

From (\ref{eq24}) and (\ref{eq25}), one can infer further that the upper bound of $\sum\limits_{k = 1}^n {\mu (k)}$, with $u = 0$ and $\sigma=\sqrt{6/\pi^2}$ for $\mu(k)$, is,

\begin{equation}\label{eq26}
\sum\limits_{k = 1}^n {\mu (k) = M(n) \le \sigma {K_{\alpha /2}}} \sqrt n  = \sqrt {6/{\pi ^2}} {K_{\alpha /2}}\sqrt n 
\end{equation}

The inequality (\ref{eq26}) holds with a probability of $1 - \alpha$.

In fact, the inequality of Mertens type is only a special case of (\ref{eq26}) with a fixed probability of $1 - \alpha$.

\section{Discussions}

\subsection{The calculation of $\mu (n) $ with the recursive relation}

In theory, we can calculate $\mu (n)$ for any large $n$ recursively with the recursive relation obtained in section \ref{sec2}. However, in order to calculate $\mu (n)$, we need to know $\mu (1), \mu (2), \mu (3),
\ldots \ldots \mu (n - 1)$. Usually $\mu (1), \mu (2), \mu (3), \ldots \ldots
\mu (n - 1)$ are stored in an array which demands much larger amount of
computer memory if $n$ is large. In this paper, we only calculate the values
of $\mu (n)$ from \ $\mu (1)$ to \ $\mu (2 \times 10^7)$ because of the memory
limitation of our desktop computer and computing time. To obtain more
numerical results of $\mu (n)$ with large $n$, both the faster and/or
optimization algorithm for the recursive relation here and better hardware
platform are required. It is a probable way by which the calculations are
divided into blocks and are computed with GPU, or quantum computer will be used in the future.

\subsection{The independent randomness of $\mu(n)$}

Based on the numerical consistency between empirical statistical quantities for only $2\times 10^7$ $\mu(n)$ and those from number theory (eg., $N(\mbox{Statistical})$) and $N(\mbox{Theorical})$, $p_{\tiny e}$ and $p_{\tiny t}$), we use classic statistical method to study $\mu (n)$ regardless of the strict validity of the independent randomness of $\mu(n)$. In the respect of the independent randomness of $\mu(n)$, there are some discussions (eg., Sarnak, 2012). Although $\mu(n)$ is deterministic from the recursive relation in section \ref{sec2}, it is visually random and independent. Figure {\ref{fig6} shows that the variation of ${\sum\limits_{k = 1}^n {\mu (k)}}/{n}$ with $n$ \ when $n =500000$. It can be observed that $\mu(n)$ has some properties of the independent random variable.

\begin{figure}[htb]
\setlength{\belowcaptionskip}{0pt}
\centering
\begin{overpic}[scale=0.44]{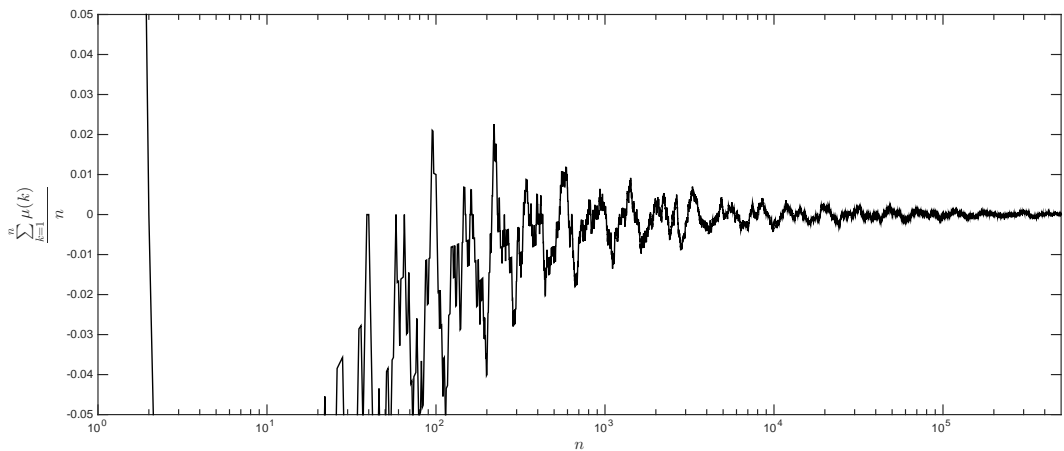}
\end{overpic}
\renewcommand{\figurename}{Fig.}
\caption{Variation of $\frac{\sum\limits_{k = 1}^n {\mu (k)}}{n}$ with $n$ \ when $n =500000$.}
\label{fig6}
\end{figure}

The above can be viewed restrictedly from another definition of $\mu(n)$,

\begin{equation}\label{eq31}
\mu(n)  = \left\{ {\begin{array}{*{20}{c}}
0\\
{{{( - 1)}^{\omega(n)}}}
\end{array}\begin{array}{*{20}{l}}
{{\mbox{if \ }} n {\mbox{\ is non-squarefree\hspace{0em}}}}\\
{{\mbox{if\  }} n {\mbox{\  is squarefree}}}
\end{array}} \right.
\end{equation}

where $\omega(n)$ is the number of distinct prime factors.

According to Erd$\ddot{\mathrm{o}}$s-Kac Theorem, $\omega(n)$ is independent and random when $n$ is large, so may be $\mu(n)$ for $\mu(n)=( - 1)^{\omega(n)}=\exp [i z \omega (n)]\vert_{_{z=\pi}}$ when $n$ is squarefree and if we take $n$ in random order.

Furthermore, $\mu(n)$ is either $(-1)^{\omega(n)}$ or $0$ because of its definition of (\ref{eq31}), so we have $M(n)=M(n')=\sum\limits_{k = 1}^{n'} {\mu (k)}=\sum\limits_{k = 1}^{n'} {(-1)^{\omega(k)}}$ ($n'$ is squarefree). The randomness of $M(n')$ should be stronger than $\mu(n)$ if we take $n'$ in random order. Numerical results of Good and Churchhouse (1968) show $M(n)$ in blocks of length $N$ has asymptotically a normal distribution with mean zero and variance of $6N/\pi^2$ (where $N$ is large). These numerical results can be rechecked as follows. From the rules in mathematical statistics, we know that the observed values of a discrete random variable $X$ ($X={x_1,x_2,\ldots,x_n}$) lie  in the following interval with probability $p> 1-\alpha$, for real number $\alpha$ with $0 < \alpha < 1$, 

\begin{equation}\label{eq32}
\left[ \overline{X}- \frac{\Delta X }{\sqrt{\alpha}},\overline{X}+ \frac{\Delta X }{\sqrt{\alpha}} \right]
\end{equation}

where $\overline{X}=\sum\limits_{k = 1}^n {x_kp_k}$, $(\Delta X)^2=\sum\limits_{k = 1}^n {(x_k-\overline{X})^2p_k}$, $p_k:=P(x=x_k)$.

Obviously, for $M(n)=M(n')=\sum\limits_{k = 1}^{n'} {(-1)^{\omega(k)}}$, $p_k=6/\pi^2$ according to subsection (\ref{subsec3_1}), and further one have,

$$\overline{X}=\frac{6}{\pi^2}\sum\limits_{k = 1}^{n'_{\mbox{\tiny even}}} {1}+\frac{6}{\pi^2}\sum\limits_{l = 1}^{n'_{\mbox{\tiny odd}}} {-1}=0$$
$$(\Delta X)^2=\frac{6}{\pi^2}\sum\limits_{k = 1}^{n'_{\mbox{\tiny even}}} {1}+\frac{6}{\pi^2}\sum\limits_{l = 1}^{n'_{\mbox{\tiny odd}}} {1}=\frac{6n'}{\pi^2}$$

where $n'_{\mbox{\tiny even}}$ is squarefree with even number of distinct prime factors, $n'_{\mbox{\tiny odd}}$  with odd number of distinct prime factors.

Therefore, from (\ref{eq32}), we can obtain an upper bound for $M(n)$ similar to (\ref{eq26}) as follows, with probability $p> 1-\alpha$,

\begin{equation}\label{eq33}
\sum\limits_{k = 1}^n {\mu (k)} = M(n)=\sum\limits_{k = 1}^{n'} {(-1)^{\omega(k)}} \le  \frac{\sqrt {6/{\pi ^2}}}{\sqrt{\alpha}} \sqrt{n'} \le  \frac{\sqrt {6/{\pi ^2}}}{\sqrt{\alpha}} \sqrt n
\end{equation}

If $\alpha$ takes $\frac{6}{\pi^2}$, then $M(n)\le\sqrt{n}$ with a probability $p>1-6/\pi^2\approx 0.3920$, which means 
that Mertens conjecture is not true.

On the other hand,  we can check whether $\mu(n)$ is periodic or not by estimating its power spectral density (PSD). We calculate the PSD for $\mu(n)$ series from $\mu(1)$ to $\mu(2\times 10^7)$ by taking $n$ as time. The results are shown in Figure \ref{fig4}.  It can be found that $\mu(n)$ has no apparent periodicity because the PSD of $\mu(n)$ have no distinguished peak(s).

\begin{figure}[htb]
\setlength{\belowcaptionskip}{0pt}
\centering
\begin{overpic}[scale=0.5]{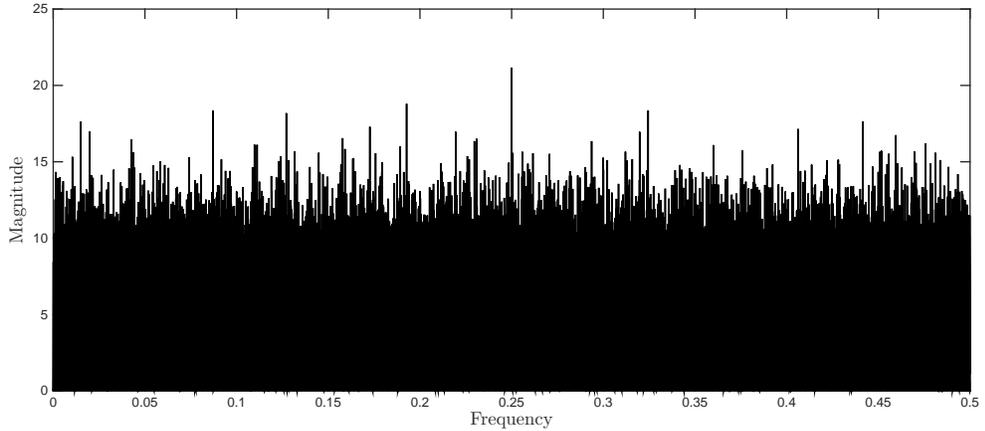}
\end{overpic}
\renewcommand{\figurename}{Fig.}
\caption{The power spectral density (PSD) for $\mu(n)$ series from $\mu(1)$ to $\mu(2\times 10^7)$.}
\label{fig4}
\end{figure}

Although the independent randomness of $\mu(n)$ is a problem unsolved so far, we can analyse $\mu (n)$ by way of statistics because $\mu (n)$ has a complicate and non-periodic distribution, as those statistical approaches applied to chaos. 

\subsection{The probability of ${\sum\limits_{k = 1}^n {\left| {\mu (k)} \right|} } > C n$}

Similarly, we can calculate the probability of ${\sum\limits_{k = 1}^n {\left| {\mu (k)} \right|} } > C n$. The distribution law for $\vert \mu (n) \vert $can be obtained as shown in Table \ref{tb2}.

\begin{table}[htdp]
\caption{The distribution rule for $\vert \mu(k)\vert$}
\begin{center}
\begin{tabular}{c c c }
\hline
$\vert \mu (k)\vert$ &$0$ &$1$\\
\hline
$p_{t_{_k}}$&$1-6/\pi^2$&$6/\pi^2$\\
\hline
\end{tabular}
\end{center}
\label{tb2}
\end{table}%

And $E (\vert \mu (k) \vert) = 6 / \pi^2$ and $D (\vert \mu (k) \vert) = 6 / \pi^2 (1 - 6 / \pi^2)$.

According to central limit theorem,  for any $x$, we have,

\begin{equation}\label{eq27}
\begin{array}{*{20}{c}}
{\mathop {\lim }\limits_{n \to \infty } P\left\{ {\frac{{\sum\limits_{k = 1}^n {\vert \mu (k)\vert - E(\sum\limits_{k = 1}^n {\vert \mu (k)\vert} )} }}{{\sqrt {D(\sum\limits_{k = 1}^n {\vert \mu (k)\vert} )} }} \le x} \right\}}\\
{}\\
{}
\end{array}
\begin{array}{*{20}{l}}
{ = \mathop {\lim }\limits_{n \to \infty } P\left\{ {\frac{{\sum\limits_{k = 1}^n {\left| {\mu (k)} \right|}  - 6n/{\pi ^2}}}{{\sqrt {n6/{\pi ^2}(1 - 6/{\pi ^2})} }} \le x} \right\}}\\
{{\rm{ }} = \int\limits_{ - \infty }^x {\frac{1}{{\sqrt {2\pi } }}\exp ( - \frac{{{t^2}}}{2}){\rm{d}}t} }\\
{{\rm{}} = \Phi (x)}
\end{array}
\end{equation}

And (\ref{eq27}) means,

\begin{equation}\label{eq28}
\frac{{\sum\limits_{k = 1}^n {\left| {\mu (k)} \right|}  - 6n/{\pi ^2}}}{{\sqrt {n6/{\pi ^2}(1 - 6/{\pi ^2})} }} \sim N(0,1)
\end{equation}

\begin{figure}[htb]
\setlength{\belowcaptionskip}{0pt}
\centering
\begin{overpic}[scale=0.6]{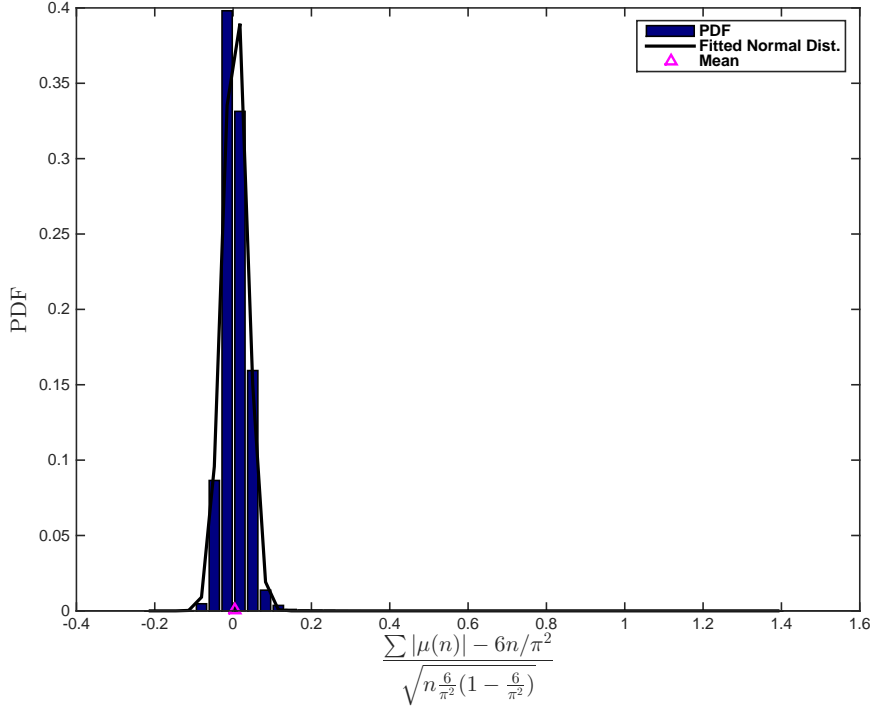}
\end{overpic}
\renewcommand{\figurename}{Fig.}
\caption{The probability density function for $\frac{\sum^n_{k = 1} \vert \mu (k) \vert - 6 n / \pi^2}{\sqrt{6 n (1 - 6 / \pi^2) / \pi^2}}$ \ when $n =500000$.}
\label{fig5}
\end{figure}

Figure \ref{fig5} shows the probability density function for $\frac{\sum^n_{k = 1} \vert \mu (k) \vert - 6 n / \pi^2}{\sqrt{6 n (1 - 6 / \pi^2) / \pi^2}}$ \ when $n =500000$. It can be seen that the distribution of (\ref{eq28}) is reasonable.

With equality (\ref{eq27}), the probability of ${\sum\limits_{k = 1}^n {\left| {\mu (k)} \right|} } > C n$ (where $C$ is a constant) can be obtained. Clearly,

\begin{equation}\label{eq29}
\begin{array}{*{20}{c}}
{P\left\{ {\sum\limits_{k = 1}^n {\left| {\mu (k)} \right|}  > Cn} \right\}}\\
{}\\
{}\\
{}\\
{}\\
{}
\end{array}
\begin{array}{*{20}{l}}
{ = P\left\{ {\frac{{\sum\limits_{k = 1}^n {\left| {\mu (k)} \right|}  - 6n/{\pi ^2}}}{{\sqrt {n6/{\pi ^2}(1 - 6/{\pi ^2})} }} > \frac{{Cn - 6n/{\pi ^2}}}{{\sqrt {n6/{\pi ^2}(1 - 6/{\pi ^2})} }}} \right\}{\rm{    }}}\\
{{\rm{}} = 1 - P\left\{ {\frac{{\sum\limits_{k = 1}^n {\left| {\mu (k)} \right|}  - 6n/{\pi ^2}}}{{\sqrt {n6/{\pi ^2}(1 - 6/{\pi ^2})} }} \le \frac{{Cn - 6n/{\pi ^2}}}{{\sqrt {n6/{\pi ^2}(1 - 6/{\pi ^2})} }}} \right\}}\\
{{\rm{}} = 1 - \Phi\left\{ {\frac{{(C - 6/\pi^2)n}}{{\sqrt {n6/{\pi ^2}(1 - 6/{\pi ^2})} }}} \right\}{\rm{}}}
\end{array}
\end{equation}

With (\ref{eq29}), when $n$ is large

\begin{equation}\label{eq30}
P\left\{ {\sum\limits_{k = 1}^n {\left| {\mu (k)} \right|}  > Cn} \right\} = \left\{ {\begin{array}{*{20}{c}}
0&{C > 6/{\pi ^2}}\\
{1/2}&{C = 6/{\pi ^2}}\\
1&{C < 6/{\pi ^2}}
\end{array}} \right.
\end{equation}

\section{Conclusions}

Based on the results and discussion above, some conclusions can be drawn as
follows,

(1) An elementary definite recursive relation for M$\ddot{\mathrm{o}}$bius function is
introduced by two simple ways. One is from M$\ddot{\mathrm{o}}$bius transform, and the
other is from Redheffer Matrix. With this recursive relation, $\mu (n)$ can be
calculated without directly knowing the factorization of $n$, in which the
most complex operation is only the Mod.

(2) With this relation, $\mu (1) \sim \mu (2 \times 10^7) $ are calculated
recursively. Based on these $2\times 10^7$ samples, we calculate the frequencies and 
empirical probabilities for  $\mu(n)$ of taking $-1,0,1$, so does for $\vert\mu\vert=1$. 
And then compare them with those in number theory. It can
be found these two kinds of frequencies and probabilities are numerically consistent.

(3) Based on these numerical results, we take $\mu (n)$ as an independent
random sequence although it has a deterministic recursive rule. The expectation
and variance of the $\mu (k)$ are $E (\mu (k)) = 0$ and $D (\mu (k)) = 6 /
\pi^2$, respectively.

(4) We show that the Mertens conjecture, even any conjecture of the Mertens
type, is false in a probability sense, and present an upper bound for cumulative
sums of $\mu (n)$ as$\overset{n}{\underset{k = 1}{\sum}}^{} \mu (k) \leqslant
\sqrt{6 / \pi^2} K_{\alpha / 2} \sqrt{n}$ with a probability of $1 - \alpha$.

\vspace{20em}

{\Large\bf  Acknowledges}

We thank Russ Woodroofe very much for pointing out that the recursive relation (\ref{eq11}) is a special case for (poset) M$\ddot{\mathrm{o}}$bius function in Incidence Algebra.

\vspace{20em}

\def\thebibliography#1{
{\Large\bf  References}\list
 {}{\setlength\labelwidth{1.4em}\leftmargin\labelwidth
 \setlength\parsep{0pt}\setlength\itemsep{.3\baselineskip}
 \setlength{\itemindent}{-\leftmargin}
 \usecounter{enumi}}
 \def\newblock{\hskip .11em plus .33em minus -.07em}
 \sloppy
 \sfcode`\.=1000\relax}
\let\endthebibliography=\endlist

\clearpage
{\Large\bf  Appendix}

\begin{longtable}{lllllll}
\caption{The comparison of empirical probability and frequency with those from number theory when $\mu(n)=-1/1$}\\
\hline
    $N$     & $N_{-1}$     & $N_{-1}^{\mbox{\tiny T}}$     & $p_{_{\tiny{e}}}$ ($0.3039635509$) & $N_1$ & $N_1^{\mbox{\tiny T}}$     & $p_{_{\tiny{e}}}$($0.3039635509$) \\
    \hline
    \endfirsthead
   \caption[]{The comparison of empirical probability and frequency with those from number theory when $\mu(n)=-1/1$ (continued)}\\
    \hline
      $N$     & $N_{-1}$     & $N_{-1}^{\mbox{\tiny T}}$     & $p_{_{\tiny{e}}}$ ($0.3039635509$) & $N_1$ & $N_1^{\mbox{\tiny T}}$     & $p_{_{\tiny{e}}}$($0.3039635509$) \\
      \hline
      \endhead
      \hline \endfoot
     100000 & 30421 & 30396.4  & 0.3042100000  & 30373 & 30396.4 0 & 0.3037300000  \\
    200000 & 60791 & 60792.7  & 0.3039550000  & 60790 & 60792.7 0 & 0.3039500000  \\
    300000 & 91079 & 91189.1  & 0.3035966667  & 91299 & 91189.1 0 & 0.3043300000  \\
    400000 & 121577 & 121585.4  & 0.3039425000  & 121588 & 121585.4  & 0.3039700000  \\
    500000 & 151982 & 151981.8  & 0.3039640000  & 151976 & 151981.8  & 0.3039520000  \\
    600000 & 182492 & 182378.1  & 0.3041533333  & 182262 & 182378.1  & 0.3037700000  \\
    700000 & 212666 & 212774.5  & 0.3038085714  & 212892 & 212774.5  & 0.3041314286  \\
    800000 & 243181 & 243170.8  & 0.3039762500  & 243161 & 243170.8  & 0.3039512500  \\
    900000 & 273678 & 273567.2  & 0.3040866667  & 273453 & 273567.2  & 0.3038366667  \\
    1000000 & 303857 & 303963.6  & 0.3038570000  & 304069 & 303963.6  & 0.3040690000  \\
    1100000 & 334263 & 334359.9  & 0.3038754545  & 334464 & 334359.9  & 0.3040581818  \\
    1200000 & 364832 & 364756.3  & 0.3040266667  & 364677 & 364756.3  & 0.3038975000  \\
    1300000 & 395192 & 395152.6  & 0.3039938462  & 395111 & 395152.6  & 0.3039315385  \\
    1400000 & 425669 & 425549.0  & 0.3040492857  & 425422 & 425549.0  & 0.3038728571  \\
    1500000 & 456092 & 455945.3  & 0.3040613333  & 455799 & 455945.3  & 0.3038660000  \\
    1600000 & 486262 & 486341.7  & 0.3039137500  & 486430 & 486341.7  & 0.3040187500  \\
    1700000 & 516688 & 516738.0  & 0.3039341176  & 516792 & 516738.0  & 0.3039952941  \\
    1800000 & 546936 & 547134.4  & 0.3038533333  & 547340 & 547134.4  & 0.3040777778  \\
    1900000 & 577533 & 577530.7  & 0.3039647368  & 577544 & 577530.7  & 0.3039705263  \\
    2000000 & 608062 & 607927.1  & 0.3040310000  & 607815 & 607927.1  & 0.3039075000  \\
    2100000 & 638508 & 638323.5  & 0.3040514286  & 638142 & 638323.5  & 0.3038771429  \\
    2200000 & 668846 & 668719.8  & 0.3040209091  & 668600 & 668719.8  & 0.3039090909  \\
    2300000 & 699203 & 699116.2  & 0.3040013043  & 699016 & 699116.2  & 0.3039200000  \\
    2400000 & 729384 & 729512.5  & 0.3039100000  & 729638 & 729512.5  & 0.3040158333  \\
    2500000 & 759726 & 759908.9  & 0.3038904000  & 760088 & 759908.9  & 0.3040352000  \\
    2600000 & 790230 & 790305.2  & 0.3039346154  & 790377 & 790305.2  & 0.3039911538  \\
    2700000 & 820674 & 820701.6  & 0.3039533333  & 820722 & 820701.6  & 0.3039711111  \\
    2800000 & 850937 & 851097.9  & 0.3039060714  & 851251 & 851097.9  & 0.3040182143  \\
    2900000 & 881525 & 881494.3  & 0.3039741379  & 881457 & 881494.3  & 0.3039506897  \\
    3000000 & 911833 & 911890.7  & 0.3039443333  & 911940 & 911890.7  & 0.3039800000  \\
    3100000 & 942195 & 942287.0  & 0.3039338710  & 942372 & 942287.0  & 0.3039909677  \\
    3200000 & 972925 & 972683.4  & 0.3040390625  & 972438 & 972683.4  & 0.3038868750  \\
    3300000 & 1003355 & 1003079.7  & 0.3040469697  & 1002803 & 1003079.7  & 0.3038796970  \\
    3400000 & 1033623 & 1033476.1  & 0.3040067647  & 1033331 & 1033476.1  & 0.3039208824  \\
    3500000 & 1063947 & 1063872.4  & 0.3039848571  & 1063809 & 1063872.4  & 0.3039454286  \\
    3600000 & 1094166 & 1094268.8  & 0.3039350000  & 1094378 & 1094268.8  & 0.3039938889  \\
    3700000 & 1124662 & 1124665.1  & 0.3039627027  & 1124661 & 1124665.1  & 0.3039624324  \\
    3800000 & 1154989 & 1155061.5  & 0.3039444737  & 1155144 & 1155061.5  & 0.3039852632  \\
    3900000 & 1185390 & 1185457.8  & 0.3039461538  & 1185540 & 1185457.8  & 0.3039846154  \\
    4000000 & 1215772 & 1215854.2  & 0.3039430000  & 1215964 & 1215854.2  & 0.3039910000  \\
    4100000 & 1246259 & 1246250.6  & 0.3039656098  & 1246258 & 1246250.6  & 0.3039653659  \\
    4200000 & 1276499 & 1276646.9  & 0.3039283333  & 1276799 & 1276646.9  & 0.3039997619  \\
    4300000 & 1306851 & 1307043.3  & 0.3039188372  & 1307252 & 1307043.3  & 0.3040120930  \\
    4400000 & 1337169 & 1337439.6  & 0.3039020455  & 1337720 & 1337439.6  & 0.3040272727  \\
    4500000 & 1367757 & 1367836.0  & 0.3039460000  & 1367930 & 1367836.0  & 0.3039844444  \\
    4600000 & 1398114 & 1398232.3  & 0.3039378261  & 1398354 & 1398232.3  & 0.3039900000  \\
    4700000 & 1428655 & 1428628.7  & 0.3039691489  & 1428604 & 1428628.7  & 0.3039582979  \\
    4800000 & 1459035 & 1459025.0  & 0.3039656250  & 1459025 & 1459025.0  & 0.3039635417  \\
    4900000 & 1489629 & 1489421.4  & 0.3040059184  & 1489223 & 1489421.4  & 0.3039230612  \\
    5000000 & 1520171 & 1519817.8  & 0.3040342000  & 1519462 & 1519817.8  & 0.3038924000  \\
    5100000 & 1550492 & 1550214.1  & 0.3040180392  & 1549940 & 1550214.1  & 0.3039098039  \\
    5200000 & 1580813 & 1580610.5  & 0.3040025000  & 1580407 & 1580610.5  & 0.3039244231  \\
    5300000 & 1611343 & 1611006.8  & 0.3040269811  & 1610658 & 1611006.8  & 0.3038977358  \\
    5400000 & 1641700 & 1641403.2  & 0.3040185185  & 1641101 & 1641403.2  & 0.3039075926  \\
  5500000  & 1672051  &1671799.5  &0.3040092727    &1671538    &1671799.5 &0.3039160000\\
5600000  & 1702293  &1702195.9  &0.3039808929    &1702098    &1702195.9 &0.3039460714\\
5700000  & 1732416  &1732592.2  &0.3039326316    &1732764    &1732592.2 &0.3039936842\\
5800000  & 1762756  &1762988.6  &0.3039234483    &1763212    &1762988.6 &0.3040020690\\
5900000  & 1793193  &1793385.0  &0.3039310169    &1793570    &1793385.0 &0.3039949153\\
6000000  & 1823650  &1823781.3  &0.3039416667    &1823907    &1823781.3 &0.3039845000\\
6100000  & 1854028  &1854177.7  &0.3039390164    &1854325    &1854177.7 &0.3039877049\\
6200000  & 1884474  &1884574.0  &0.3039474194    &1884666    &1884574.0 &0.3039783871\\
6300000  & 1914862  &1914970.4  &0.3039463492    &1915064    &1914970.4 &0.3039784127\\
6400000  & 1945014  &1945366.7  &0.3039084375    &1945715    &1945366.7 &0.3040179687\\
6500000  & 1975328  &1975763.1  &0.3038966154    &1976195    &1975763.1 &0.3040300000\\
6600000  & 2005813  &2006159.4  &0.3039110606    &2006509    &2006159.4 &0.3040165152\\
6700000  & 2036280  &2036555.8  &0.3039223881    &2036837    &2036555.8 &0.3040055224\\
6800000  & 2066717  &2066952.1  &0.3039289706    &2067210    &2066952.1 &0.3040014706\\
6900000  & 2097200  &2097348.5  &0.3039420290    &2097505    &2097348.5 &0.3039862319\\
7000000  & 2127844  &2127744.9  &0.3039777143    &2127660    &2127744.9 &0.3039514286\\
7100000  & 2158600  &2158141.2  &0.3040281690    &2157693    &2158141.2 &0.3039004225\\
7200000  & 2188861  &2188537.6  &0.3040084722    &2188217    &2188537.6 &0.3039190278\\
7300000  & 2219265  &2218933.9  &0.3040089041    &2218604    &2218933.9 &0.3039183562\\
7400000  & 2249572  &2249330.3  &0.3039962162    &2249081    &2249330.3 &0.3039298649\\
7500000  & 2279804  &2279726.6  &0.3039738667    &2279632    &2279726.6 &0.3039509333\\
7600000  & 2310212  &2310123.0  &0.3039752632    &2310014    &2310123.0 &0.3039492105\\
7700000  & 2340655  &2340519.3  &0.3039811688    &2340371    &2340519.3 &0.3039442857\\
7800000  & 2370879  &2370915.7  &0.3039588462    &2370945    &2370915.7 &0.3039673077\\
7900000  & 2401416  &2401312.1  &0.3039767089    &2401202    &2401312.1 &0.3039496203\\
8000000  & 2431796  &2431708.4  &0.3039745000    &2431607    &2431708.4 &0.3039508750\\
8100000  & 2462041  &2462104.8  &0.3039556790    &2462157    &2462104.8 &0.3039700000\\
8200000  & 2492428  &2492501.1  &0.3039546341    &2492563    &2492501.1 &0.3039710976\\
8300000  & 2522880  &2522897.5  &0.3039614458    &2522900    &2522897.5 &0.3039638554\\
8400000  & 2553329  &2553293.8  &0.3039677381    &2553250    &2553293.8 &0.3039583333\\
8500000  & 2583583  &2583690.2  &0.3039509412    &2583782    &2583690.2 &0.3039743529\\
8600000  & 2614150  &2614086.5  &0.3039709302    &2614001    &2614086.5 &0.3039536047\\
8700000  & 2644785  &2644482.9  &0.3039982759    &2644163    &2644482.9 &0.3039267816\\
8800000  & 2675324  &2674879.2  &0.3040140909    &2674412    &2674879.2 &0.3039104545\\
8900000  & 2705512  &2705275.6  &0.3039901124    &2705036    &2705275.6 &0.3039366292\\
9000000  & 2735841  &2735672.0  &0.3039823333    &2735501    &2735672.0 &0.3039445556\\
9100000  & 2766275  &2766068.3  &0.3039862637    &2765847    &2766068.3 &0.3039392308\\
9200000  & 2796640  &2796464.7  &0.3039826087    &2796280    &2796464.7 &0.3039434783\\
9300000  & 2827009  &2826861.0  &0.3039794624    &2826712    &2826861.0 &0.3039475269\\
9400000  & 2857367  &2857257.4  &0.3039752128    &2857156    &2857257.4 &0.3039527660\\
9500000  & 2887605  &2887653.7  &0.3039584211    &2887723    &2887653.7 &0.3039708421\\
9600000  & 2918034  &2918050.1  &0.3039618750    &2918074    &2918050.1 &0.3039660417\\
9700000  & 2948319  &2948446.4  &0.3039504124    &2948598    &2948446.4 &0.3039791753\\
9800000  & 2978563  &2978842.8  &0.3039350000    &2979144    &2978842.8 &0.3039942857\\
9900000  & 3008889  &3009239.2  &0.3039281818    &3009605    &3009239.2 &0.3040005051\\
10000000 & 3039127  &3039635.5  &0.3039127000    &3040164    &3039635.5 &0.3040164000\\
10100000 & 3069687  &3070031.9  &0.3039294059    &3070405    &3070031.9 &0.3040004950\\
10200000 & 3099872  &3100428.2  &0.3039090196    &3101015    &3100428.2 &0.3040210784\\
10300000 & 3130336  &3130824.6  &0.3039161165    &3131328    &3130824.6 &0.3040124272\\
10400000 & 3160738  &3161220.9  &0.3039171154    &3161719    &3161220.9 &0.3040114423\\
10500000 & 3191206  &3191617.3  &0.3039243810    &3192048    &3191617.3 &0.3040045714\\
10600000 & 3221686  &3222013.6  &0.3039326415    &3222348    &3222013.6 &0.3039950943\\
10700000 & 3252202  &3252410.0  &0.3039441121    &3252627    &3252410.0 &0.3039838318\\
10800000 & 3282578  &3282806.4  &0.3039424074    &3283049    &3282806.4 &0.3039860185\\
10900000 & 3312999  &3313202.7  &0.3039448624    &3313429    &3313202.7 &0.3039843119\\
11000000 & 3343506  &3343599.1  &0.3039550909    &3343712    &3343599.1 &0.3039738182\\
11100000 & 3374011  &3373995.4  &0.3039649550    &3373989    &3373995.4 &0.3039629730\\
11200000 & 3404419  &3404391.8  &0.3039659821    &3404376    &3404391.8 &0.3039621429\\
11300000 & 3434879  &3434788.1  &0.3039715929    &3434711    &3434788.1 &0.3039567257\\
11400000 & 3465185  &3465184.5  &0.3039635965    &3465208    &3465184.5 &0.3039656140\\
11500000 & 3495628  &3495580.8  &0.3039676522    &3495538    &3495580.8 &0.3039598261\\
11600000 & 3526180  &3525977.2  &0.3039810345    &3525787    &3525977.2 &0.3039471552\\
11700000 & 3556685  &3556373.5  &0.3039901709    &3556070    &3556373.5 &0.3039376068\\
11800000 & 3587067  &3586769.9  &0.3039887288    &3586490    &3586769.9 &0.3039398305\\
11900000 & 3617377  &3617166.3  &0.3039812605    &3616976    &3617166.3 &0.3039475630\\
12000000 & 3647897  &3647562.6  &0.3039914167    &3647243    &3647562.6 &0.3039369167\\
12100000 & 3678249  &3677959.0  &0.3039875207    &3677691    &3677959.0 &0.3039414050\\
12200000 & 3708727  &3708355.3  &0.3039940164    &3707994    &3708355.3 &0.3039339344\\
12300000 & 3738956  &3738751.7  &0.3039801626    &3738560    &3738751.7 &0.3039479675\\
12400000 & 3769655  &3769148.0  &0.3040044355    &3768649    &3769148.0 &0.3039233065\\
12500000 & 3800038  &3799544.4  &0.3040030400    &3799058    &3799544.4 &0.3039246400\\
12600000 & 3830284  &3829940.7  &0.3039907937    &3829590    &3829940.7 &0.3039357143\\
12700000 & 3860722  &3860337.1  &0.3039938583    &3859952    &3860337.1 &0.3039332283\\
12800000 & 3891316  &3890733.5  &0.3040090625    &3890157    &3890733.5 &0.3039185156\\
12900000 & 3921803  &3921129.8  &0.3040157364    &3920466    &3921129.8 &0.3039120930\\
13000000 & 3952009  &3951526.2  &0.3040006923    &3951058    &3951526.2 &0.3039275385\\
13100000 & 3982303  &3981922.5  &0.3039925954    &3981547    &3981922.5 &0.3039348855\\
13200000 & 4012618  &4012318.9  &0.3039862121    &4012019    &4012318.9 &0.3039408333\\
13300000 & 4042856  &4042715.2  &0.3039741353    &4042567    &4042715.2 &0.3039524060\\
13400000 & 4073262  &4073111.6  &0.3039747761    &4072950    &4073111.6 &0.3039514925\\
13500000 & 4103423  &4103507.9  &0.3039572593    &4103593    &4103507.9 &0.3039698519\\
13600000 & 4133881  &4133904.3  &0.3039618382    &4133938    &4133904.3 &0.3039660294\\
13700000 & 4164370  &4164300.6  &0.3039686131    &4164238    &4164300.6 &0.3039589781\\
13800000 & 4194724  &4194697.0  &0.3039655072    &4194671    &4194697.0 &0.3039616667\\
13900000 & 4225100  &4225093.4  &0.3039640288    &4225097    &4225093.4 &0.3039638129\\
14000000 & 4255463  &4255489.7  &0.3039616429    &4255509    &4255489.7 &0.3039649286\\
14100000 & 4285549  &4285886.1  &0.3039396454    &4286231    &4285886.1 &0.3039880142\\
14200000 & 4315867  &4316282.4  &0.3039342958    &4316712    &4316282.4 &0.3039938028\\
14300000 & 4346401  &4346678.8  &0.3039441259    &4346971    &4346678.8 &0.3039839860\\
14400000 & 4377024  &4377075.1  &0.3039600000    &4377125    &4377075.1 &0.3039670139\\
14500000 & 4407257  &4407471.5  &0.3039487586    &4407682    &4407471.5 &0.3039780690\\
14600000 & 4437734  &4437867.8  &0.3039543836    &4438008    &4437867.8 &0.3039731507\\
14700000 & 4467979  &4468264.2  &0.3039441497    &4468549    &4468264.2 &0.3039829252\\
14800000 & 4498441  &4498660.6  &0.3039487162    &4498892    &4498660.6 &0.3039791892\\
14900000 & 4528728  &4529056.9  &0.3039414765    &4529378    &4529056.9 &0.3039851007\\
15000000 & 4559112  &4559453.3  &0.3039408000    &4559777    &4559453.3 &0.3039851333\\
15100000 & 4589543  &4589849.6  &0.3039432450    &4590147    &4589849.6 &0.3039832450\\
15200000 & 4619957  &4620246.0  &0.3039445395    &4620528    &4620246.0 &0.3039821053\\
15300000 & 4650383  &4650642.3  &0.3039466013    &4650883    &4650642.3 &0.3039792810\\
15400000 & 4680596  &4681038.7  &0.3039348052    &4681468    &4681038.7 &0.3039914286\\
15500000 & 4711047  &4711435.0  &0.3039385161    &4711798    &4711435.0 &0.3039869677\\
15600000 & 4741674  &4741831.4  &0.3039534615    &4741957    &4741831.4 &0.3039716026\\
15700000 & 4771945  &4772227.7  &0.3039455414    &4772481    &4772227.7 &0.3039796815\\
15800000 & 4802182  &4802624.1  &0.3039355696    &4803037    &4802624.1 &0.3039896835\\
15900000 & 4832633  &4833020.5  &0.3039391824    &4833375    &4833020.5 &0.3039858491\\
16000000 & 4863066  &4863416.8  &0.3039416250    &4863730    &4863416.8 &0.3039831250\\
16100000 & 4893413  &4893813.2  &0.3039386957    &4894182    &4893813.2 &0.3039864596\\
16200000 & 4923786  &4924209.5  &0.3039374074    &4924602    &4924209.5 &0.3039877778\\
16300000 & 4954368  &4954605.9  &0.3039489571    &4954828    &4954605.9 &0.3039771779\\
16400000 & 4984486  &4985002.2  &0.3039320732    &4985493    &4985002.2 &0.3039934756\\
16500000 & 5014984  &5015398.6  &0.3039384242    &5015787    &5015398.6 &0.3039870909\\
16600000 & 5045530  &5045794.9  &0.3039475904    &5046040    &5045794.9 &0.3039783133\\
16700000 & 5075930  &5076191.3  &0.3039479042    &5076427    &5076191.3 &0.3039776647\\
16800000 & 5106482  &5106587.7  &0.3039572619    &5106673    &5106587.7 &0.3039686310\\
16900000 & 5137006  &5136984.0  &0.3039648521    &5136938    &5136984.0 &0.3039608284\\
17000000 & 5167558  &5167380.4  &0.3039740000    &5167181    &5167380.4 &0.3039518235\\
17100000 & 5197978  &5197776.7  &0.3039753216    &5197540    &5197776.7 &0.3039497076\\
17200000 & 5228079  &5228173.1  &0.3039580814    &5228214    &5228173.1 &0.3039659302\\
17300000 & 5258540  &5258569.4  &0.3039618497    &5258563    &5258569.4 &0.3039631792\\
17400000 & 5289027  &5288965.8  &0.3039670690    &5288881    &5288965.8 &0.3039586782\\
17500000 & 5319647  &5319362.1  &0.3039798286    &5319054    &5319362.1 &0.3039459429\\
17600000 & 5350010  &5349758.5  &0.3039778409    &5349482    &5349758.5 &0.3039478409\\
17700000 & 5380420  &5380154.9  &0.3039785311    &5379862    &5380154.9 &0.3039470056\\
17800000 & 5410822  &5410551.2  &0.3039787640    &5410261    &5410551.2 &0.3039472472\\
17900000 & 5441148  &5440947.6  &0.3039747486    &5440722    &5440947.6 &0.3039509497\\
18000000 & 5471659  &5471343.9  &0.3039810556    &5470992    &5471343.9 &0.3039440000\\
18100000 & 5501925  &5501740.3  &0.3039737569    &5501525    &5501740.3 &0.3039516575\\
18200000 & 5532231  &5532136.6  &0.3039687363    &5532004    &5532136.6 &0.3039562637\\
18300000 & 5562820  &5562533.0  &0.3039792350    &5562211    &5562533.0 &0.3039459563\\
18400000 & 5593139  &5592929.3  &0.3039749457    &5592689    &5592929.3 &0.3039504891\\
18500000 & 5623413  &5623325.7  &0.3039682703    &5623208    &5623325.7 &0.3039571892\\
18600000 & 5653791  &5653722.0  &0.3039672581    &5653631    &5653722.0 &0.3039586559\\
18700000 & 5684183  &5684118.4  &0.3039670053    &5684045    &5684118.4 &0.3039596257\\
18800000 & 5714612  &5714514.8  &0.3039687234    &5714404    &5714514.8 &0.3039576596\\
18900000 & 5745055  &5744911.1  &0.3039711640    &5744769    &5744911.1 &0.3039560317\\
19000000 & 5775469  &5775307.5  &0.3039720526    &5775143    &5775307.5 &0.3039548947\\
19100000 & 5805915  &5805703.8  &0.3039746073    &5805503    &5805703.8 &0.3039530366\\
19200000 & 5836335  &5836100.2  &0.3039757813    &5835864    &5836100.2 &0.3039512500\\
19300000 & 5866727  &5866496.5  &0.3039754922    &5866274    &5866496.5 &0.3039520207\\
19400000 & 5897268  &5896892.9  &0.3039828866    &5896543    &5896892.9 &0.3039455155\\
19500000 & 5927598  &5927289.2  &0.3039793846    &5926989    &5927289.2 &0.3039481538\\
19600000 & 5957965  &5957685.6  &0.3039778061    &5957419    &5957685.6 &0.3039499490\\
19700000 & 5988337  &5988082.0  &0.3039764975    &5987834    &5988082.0 &0.3039509645\\
19800000 & 6018757  &6018478.3  &0.3039776263    &6018221    &6018478.3 &0.3039505556\\
19900000 & 6049168  &6048874.7  &0.3039782915    &6048608    &6048874.7 &0.3039501508\\
20000000 & 6079764  &6079271.0  &0.3039882000    &6078811    &6079271.0 &0.3039405500\\
        \hline
\caption*{\footnotesize Note:  $N$: Length of the block from $\mu(1)$ to $\mu(N)$; $N_{-1}$: Frequency of $\mu(n)=-1$ ;  $N_{-1}^{\mbox{\tiny T}}(=N\times p_{_{\tiny{t}}})$: Frequency of $\mu(n)=-1$ from number theory; $p_{_{\tiny{e}}}(=\frac{N_{-1}}{N})$: Empirical probability; The number in bracket is theoretical probability from number theory $p_{_{\tiny{t}}}$.  $N_1$ and $N_1^{\mbox{\tiny T}}$ are those for $\mu(n)=1$. }
\label{tb3}
    \end{longtable}
   
\clearpage    
\begin{longtable}{lllllll}
\caption{The comparison of empirical probability and frequency with those from number theory when $\mu(n)=0/\vert\mu (n)\vert =1$}\\
\hline
    $N$     & $N_0$     & $N_0^{\mbox{\tiny T}}$     & $p_{_{\tiny{e}}}$ ($0.3920728981$) & $N_{\vert\mu(n)\vert=1}$ & $N_{_{\vert\mu(n)\vert=1}}^{\mbox{\tiny T}}$     & $p_{_{\tiny{e}}}$(0.6079271019) \\
    \hline
    \endfirsthead
   \caption[]{The comparison of empirical probability and frequency with those from number theory when $\mu(n)=0/\vert\mu (n)\vert =1$ (continued)}\\
    \hline
      $N$     & $N_0$     & $N_0^{\mbox{\tiny T}}$     & $p_{_{\tiny{e}}}$ ($0.3920728981$) & $N_{\vert\mu(n)\vert=1}$ & $N_{\vert\mu(n)\vert=1}^{\mbox{\tiny T}}$     & $p_{_{\tiny{e}}}$($0.6079271019$) \\
      \hline
      \endhead
      \hline \endfoot
    100000 & 39206 & 39207.3  & 0.3920600000  & 60794 & 60792.7  & 0.6079400000  \\
    200000 & 78419 & 78414.6  & 0.3920950000  & 121581 & 121585.4  & 0.6079050000  \\
    300000 & 117622 & 117621.9  & 0.3920733333  & 182378 & 182378.1  & 0.6079266667  \\
    400000 & 156835 & 156829.2  & 0.3920875000  & 243165 & 243170.8  & 0.6079125000  \\
    500000 & 196042 & 196036.4  & 0.3920840000  & 303958 & 303963.6  & 0.6079160000  \\
    600000 & 235246 & 235243.7  & 0.3920766667  & 364754 & 364756.3  & 0.6079233333  \\
    700000 & 274442 & 274451.0  & 0.3920600000  & 425558 & 425549.0  & 0.6079400000  \\
    800000 & 313658 & 313658.3  & 0.3920725000  & 486342 & 486341.7  & 0.6079275000  \\
    900000 & 352869 & 352865.6  & 0.3920766667  & 547131 & 547134.4  & 0.6079233333  \\
    1000000 & 392074 & 392072.9  & 0.3920740000  & 607926 & 607927.1  & 0.6079260000  \\
    1100000 & 431273 & 431280.2  & 0.3920663636  & 668727 & 668719.8  & 0.6079336364  \\
    1200000 & 470491 & 470487.5  & 0.3920758333  & 729509 & 729512.5  & 0.6079241667  \\
    1300000 & 509697 & 509694.8  & 0.3920746154  & 790303 & 790305.2  & 0.6079253846  \\
    1400000 & 548909 & 548902.1  & 0.3920778571  & 851091 & 851097.9  & 0.6079221429  \\
    1500000 & 588109 & 588109.3  & 0.3920726667  & 911891 & 911890.7  & 0.6079273333  \\
    1600000 & 627308 & 627316.6  & 0.3920675000  & 972692 & 972683.4  & 0.6079325000  \\
    1700000 & 666520 & 666523.9  & 0.3920705882  & 1033480 & 1033476.1  & 0.6079294118  \\
    1800000 & 705724 & 705731.2  & 0.3920688889  & 1094276 & 1094268.8  & 0.6079311111  \\
    1900000 & 744923 & 744938.5  & 0.3920647368  & 1155077 & 1155061.5  & 0.6079352632  \\
    2000000 & 784123 & 784145.8  & 0.3920615000  & 1215877 & 1215854.2  & 0.6079385000  \\
    2100000 & 823350 & 823353.1  & 0.3920714286  & 1276650 & 1276646.9  & 0.6079285714  \\
    2200000 & 862554 & 862560.4  & 0.3920700000  & 1337446 & 1337439.6  & 0.6079300000  \\
    2300000 & 901781 & 901767.7  & 0.3920786957  & 1398219 & 1398232.3  & 0.6079213043  \\
    2400000 & 940978 & 940975.0  & 0.3920741667  & 1459022 & 1459025.0  & 0.6079258333  \\
    2500000 & 980186 & 980182.2  & 0.3920744000  & 1519814 & 1519817.8  & 0.6079256000  \\
    2600000 & 1019393 & 1019389.5  & 0.3920742308  & 1580607 & 1580610.5  & 0.6079257692  \\
    2700000 & 1058604 & 1058596.8  & 0.3920755556  & 1641396 & 1641403.2  & 0.6079244444  \\
    2800000 & 1097812 & 1097804.1  & 0.3920757143  & 1702188 & 1702195.9  & 0.6079242857  \\
    2900000 & 1137018 & 1137011.4  & 0.3920751724  & 1762982 & 1762988.6  & 0.6079248276  \\
    3000000 & 1176227 & 1176218.7  & 0.3920756667  & 1823773 & 1823781.3  & 0.6079243333  \\
    3100000 & 1215433 & 1215426.0  & 0.3920751613  & 1884567 & 1884574.0  & 0.6079248387  \\
    3200000 & 1254637 & 1254633.3  & 0.3920740625  & 1945363 & 1945366.7  & 0.6079259375  \\
    3300000 & 1293842 & 1293840.6  & 0.3920733333  & 2006158 & 2006159.4  & 0.6079266667  \\
    3400000 & 1333046 & 1333047.9  & 0.3920723529  & 2066954 & 2066952.1  & 0.6079276471  \\
    3500000 & 1372244 & 1372255.1  & 0.3920697143  & 2127756 & 2127744.9  & 0.6079302857  \\
    3600000 & 1411456 & 1411462.4  & 0.3920711111  & 2188544 & 2188537.6  & 0.6079288889  \\
    3700000 & 1450677 & 1450669.7  & 0.3920748649  & 2249323 & 2249330.3  & 0.6079251351  \\
    3800000 & 1489867 & 1489877.0  & 0.3920702632  & 2310133 & 2310123.0  & 0.6079297368  \\
    3900000 & 1529070 & 1529084.3  & 0.3920692308  & 2370930 & 2370915.7  & 0.6079307692  \\
    4000000 & 1568264 & 1568291.6  & 0.3920660000  & 2431736 & 2431708.4  & 0.6079340000  \\
    4100000 & 1607483 & 1607498.9  & 0.3920690244  & 2492517 & 2492501.1  & 0.6079309756  \\
    4200000 & 1646702 & 1646706.2  & 0.3920719048  & 2553298 & 2553293.8  & 0.6079280952  \\
    4300000 & 1685897 & 1685913.5  & 0.3920690698  & 2614103 & 2614086.5  & 0.6079309302  \\
    4400000 & 1725111 & 1725120.8  & 0.3920706818  & 2674889 & 2674879.2  & 0.6079293182  \\
    4500000 & 1764313 & 1764328.0  & 0.3920695556  & 2735687 & 2735672.0  & 0.6079304444  \\
    4600000 & 1803532 & 1803535.3  & 0.3920721739  & 2796468 & 2796464.7  & 0.6079278261  \\
    4700000 & 1842741 & 1842742.6  & 0.3920725532  & 2857259 & 2857257.4  & 0.6079274468  \\
    4800000 & 1881940 & 1881949.9  & 0.3920708333  & 2918060 & 2918050.1  & 0.6079291667  \\
    4900000 & 1921148 & 1921157.2  & 0.3920710204  & 2978852 & 2978842.8  & 0.6079289796  \\
    5000000 & 1960367 & 1960364.5  & 0.3920734000  & 3039633 & 3039635.5  & 0.6079266000  \\
    5100000 & 1999568 & 1999571.8  & 0.3920721569  & 3100432 & 3100428.2  & 0.6079278431  \\
    5200000 & 2038780 & 2038779.1  & 0.3920730769  & 3161220 & 3161220.9  & 0.6079269231  \\
    5300000 & 2077999 & 2077986.4  & 0.3920752830  & 3222001 & 3222013.6  & 0.6079247170  \\
    5400000 & 2117199 & 2117193.6  & 0.3920738889  & 3282801 & 3282806.4  & 0.6079261111  \\
  5500000  &  2156411    &2156400.9 &0.3920747273 & 3343589    &3343599.1  &0.6079252727 \\ 
5600000  &  2195609    &2195608.2 &0.3920730357 & 3404391    &3404391.8  &0.6079269643 \\ 
5700000  &  2234820    &2234815.5 &0.3920736842 & 3465180    &3465184.5  &0.6079263158 \\ 
5800000  &  2274032    &2274022.8 &0.3920744828 & 3525968    &3525977.2  &0.6079255172 \\ 
5900000  &  2313237    &2313230.1 &0.3920740678 & 3586763    &3586769.9  &0.6079259322 \\ 
6000000  &  2352443    &2352437.4 &0.3920738333 & 3647557    &3647562.6  &0.6079261667 \\ 
6100000  &  2391647    &2391644.7 &0.3920732787 & 3708353    &3708355.3  &0.6079267213 \\ 
6200000  &  2430860    &2430852.0 &0.3920741935 & 3769140    &3769148.0  &0.6079258065 \\ 
6300000  &  2470074    &2470059.3 &0.3920752381 & 3829926    &3829940.7  &0.6079247619 \\ 
6400000  &  2509271    &2509266.5 &0.3920735938 & 3890729    &3890733.5  &0.6079264063 \\ 
6500000  &  2548477    &2548473.8 &0.3920733846 & 3951523    &3951526.2  &0.6079266154 \\ 
6600000  &  2587678    &2587681.1 &0.3920724242 & 4012322    &4012318.9  &0.6079275758 \\ 
6700000  &  2626883    &2626888.4 &0.3920720896 & 4073117    &4073111.6  &0.6079279104 \\ 
6800000  &  2666073    &2666095.7 &0.3920695588 & 4133927    &4133904.3  &0.6079304412 \\ 
6900000  &  2705295    &2705303.0 &0.3920717391 & 4194705    &4194697.0  &0.6079282609 \\ 
7000000  &  2744496    &2744510.3 &0.3920708571 & 4255504    &4255489.7  &0.6079291429 \\ 
7100000  &  2783707    &2783717.6 &0.3920714085 & 4316293    &4316282.4  &0.6079285915 \\ 
7200000  &  2822922    &2822924.9 &0.3920725000 & 4377078    &4377075.1  &0.6079275000 \\ 
7300000  &  2862131    &2862132.2 &0.3920727397 & 4437869    &4437867.8  &0.6079272603 \\ 
7400000  &  2901347    &2901339.4 &0.3920739189 & 4498653    &4498660.6  &0.6079260811 \\ 
7500000  &  2940564    &2940546.7 &0.3920752000 & 4559436    &4559453.3  &0.6079248000 \\ 
7600000  &  2979774    &2979754.0 &0.3920755263 & 4620226    &4620246.0  &0.6079244737 \\ 
7700000  &  3018974    &3018961.3 &0.3920745455 & 4681026    &4681038.7  &0.6079254545 \\ 
7800000  &  3058176    &3058168.6 &0.3920738462 & 4741824    &4741831.4  &0.6079261538 \\ 
7900000  &  3097382    &3097375.9 &0.3920736709 & 4802618    &4802624.1  &0.6079263291 \\ 
8000000  &  3136597    &3136583.2 &0.3920746250 & 4863403    &4863416.8  &0.6079253750 \\ 
8100000  &  3175802    &3175790.5 &0.3920743210 & 4924198    &4924209.5  &0.6079256790 \\ 
8200000  &  3215009    &3214997.8 &0.3920742683 & 4984991    &4985002.2  &0.6079257317 \\ 
8300000  &  3254220    &3254205.1 &0.3920746988 & 5045780    &5045794.9  &0.6079253012 \\ 
8400000  &  3293421    &3293412.3 &0.3920739286 & 5106579    &5106587.7  &0.6079260714 \\ 
8500000  &  3332635    &3332619.6 &0.3920747059 & 5167365    &5167380.4  &0.6079252941 \\ 
8600000  &  3371849    &3371826.9 &0.3920754651 & 5228151    &5228173.1  &0.6079245349 \\ 
8700000  &  3411052    &3411034.2 &0.3920749425 & 5288948    &5288965.8  &0.6079250575 \\ 
8800000  &  3450264    &3450241.5 &0.3920754545 & 5349736    &5349758.5  &0.6079245455 \\ 
8900000  &  3489452    &3489448.8 &0.3920732584 & 5410548    &5410551.2  &0.6079267416 \\ 
9000000  &  3528658    &3528656.1 &0.3920731111 & 5471342    &5471343.9  &0.6079268889 \\ 
9100000  &  3567878    &3567863.4 &0.3920745055 & 5532122    &5532136.6  &0.6079254945 \\ 
9200000  &  3607080    &3607070.7 &0.3920739130 & 5592920    &5592929.3  &0.6079260870 \\ 
9300000  &  3646279    &3646278.0 &0.3920730108 & 5653721    &5653722.0  &0.6079269892 \\ 
9400000  &  3685477    &3685485.2 &0.3920720213 & 5714523    &5714514.8  &0.6079279787 \\ 
9500000  &  3724672    &3724692.5 &0.3920707368 & 5775328    &5775307.5  &0.6079292632 \\ 
9600000  &  3763892    &3763899.8 &0.3920720833 & 5836108    &5836100.2  &0.6079279167 \\ 
9700000  &  3803083    &3803107.1 &0.3920704124 & 5896917    &5896892.9  &0.6079295876 \\ 
9800000  &  3842293    &3842314.4 &0.3920707143 & 5957707    &5957685.6  &0.6079292857 \\ 
9900000  &  3881506    &3881521.7 &0.3920713131 & 6018494    &6018478.3  &0.6079286869 \\ 
10000000 &  3920709    &3920729.0 &0.3920709000 & 6079291    &6079271.0  &0.6079291000 \\ 
10100000 &  3959908    &3959936.3 &0.3920700990 & 6140092    &6140063.7  &0.6079299010 \\ 
10200000 &  3999113    &3999143.6 &0.3920699020 & 6200887    &6200856.4  &0.6079300980 \\ 
10300000 &  4038336    &4038350.9 &0.3920714563 & 6261664    &6261649.1  &0.6079285437 \\ 
10400000 &  4077543    &4077558.1 &0.3920714423 & 6322457    &6322441.9  &0.6079285577 \\ 
10500000 &  4116746    &4116765.4 &0.3920710476 & 6383254    &6383234.6  &0.6079289524 \\ 
10600000 &  4155966    &4155972.7 &0.3920722642 & 6444034    &6444027.3  &0.6079277358 \\ 
10700000 &  4195171    &4195180.0 &0.3920720561 & 6504829    &6504820.0  &0.6079279439 \\ 
10800000 &  4234373    &4234387.3 &0.3920715741 & 6565627    &6565612.7  &0.6079284259 \\ 
10900000 &  4273572    &4273594.6 &0.3920708257 & 6626428    &6626405.4  &0.6079291743 \\ 
11000000 &  4312782    &4312801.9 &0.3920710909 & 6687218    &6687198.1  &0.6079289091 \\ 
11100000 &  4352000    &4352009.2 &0.3920720721 & 6748000    &6747990.8  &0.6079279279 \\ 
11200000 &  4391205    &4391216.5 &0.3920718750 & 6808795    &6808783.5  &0.6079281250 \\ 
11300000 &  4430410    &4430423.7 &0.3920716814 & 6869590    &6869576.3  &0.6079283186 \\ 
11400000 &  4469607    &4469631.0 &0.3920707895 & 6930393    &6930369.0  &0.6079292105 \\ 
11500000 &  4508834    &4508838.3 &0.3920725217 & 6991166    &6991161.7  &0.6079274783 \\ 
11600000 &  4548033    &4548045.6 &0.3920718103 & 7051967    &7051954.4  &0.6079281897 \\ 
11700000 &  4587245    &4587252.9 &0.3920722222 & 7112755    &7112747.1  &0.6079277778 \\ 
11800000 &  4626443    &4626460.2 &0.3920714407 & 7173557    &7173539.8  &0.6079285593 \\ 
11900000 &  4665647    &4665667.5 &0.3920711765 & 7234353    &7234332.5  &0.6079288235 \\ 
12000000 &  4704860    &4704874.8 &0.3920716667 & 7295140    &7295125.2  &0.6079283333 \\ 
12100000 &  4744060    &4744082.1 &0.3920710744 & 7355940    &7355917.9  &0.6079289256 \\ 
12200000 &  4783279    &4783289.4 &0.3920720492 & 7416721    &7416710.6  &0.6079279508 \\ 
12300000 &  4822484    &4822496.6 &0.3920718699 & 7477516    &7477503.4  &0.6079281301 \\ 
12400000 &  4861696    &4861703.9 &0.3920722581 & 7538304    &7538296.1  &0.6079277419 \\ 
12500000 &  4900904    &4900911.2 &0.3920723200 & 7599096    &7599088.8  &0.6079276800 \\ 
12600000 &  4940126    &4940118.5 &0.3920734921 & 7659874    &7659881.5  &0.6079265079 \\ 
12700000 &  4979326    &4979325.8 &0.3920729134 & 7720674    &7720674.2  &0.6079270866 \\ 
12800000 &  5018527    &5018533.1 &0.3920724219 & 7781473    &7781466.9  &0.6079275781 \\ 
12900000 &  5057731    &5057740.4 &0.3920721705 & 7842269    &7842259.6  &0.6079278295 \\ 
13000000 &  5096933    &5096947.7 &0.3920717692 & 7903067    &7903052.3  &0.6079282308 \\ 
13100000 &  5136150    &5136155.0 &0.3920725191 & 7963850    &7963845.0  &0.6079274809 \\ 
13200000 &  5175363    &5175362.3 &0.3920729545 & 8024637    &8024637.7  &0.6079270455 \\ 
13300000 &  5214577    &5214569.5 &0.3920734586 & 8085423    &8085430.5  &0.6079265414 \\ 
13400000 &  5253788    &5253776.8 &0.3920737313 & 8146212    &8146223.2  &0.6079262687 \\ 
13500000 &  5292984    &5292984.1 &0.3920728889 & 8207016    &8207015.9  &0.6079271111 \\ 
13600000 &  5332181    &5332191.4 &0.3920721324 & 8267819    &8267808.6  &0.6079278676 \\ 
13700000 &  5371392    &5371398.7 &0.3920724088 & 8328608    &8328601.3  &0.6079275912 \\ 
13800000 &  5410605    &5410606.0 &0.3920728261 & 8389395    &8389394.0  &0.6079271739 \\ 
13900000 &  5449803    &5449813.3 &0.3920721583 & 8450197    &8450186.7  &0.6079278417 \\ 
14000000 &  5489028    &5489020.6 &0.3920734286 & 8510972    &8510979.4  &0.6079265714 \\ 
14100000 &  5528220    &5528227.9 &0.3920723404 & 8571780    &8571772.1  &0.6079276596 \\ 
14200000 &  5567421    &5567435.2 &0.3920719014 & 8632579    &8632564.8  &0.6079280986 \\ 
14300000 &  5606628    &5606642.4 &0.3920718881 & 8693372    &8693357.6  &0.6079281119 \\ 
14400000 &  5645851    &5645849.7 &0.3920729861 & 8754149    &8754150.3  &0.6079270139 \\ 
14500000 &  5685061    &5685057.0 &0.3920731724 & 8814939    &8814943.0  &0.6079268276 \\ 
14600000 &  5724258    &5724264.3 &0.3920724658 & 8875742    &8875735.7  &0.6079275342 \\ 
14700000 &  5763472    &5763471.6 &0.3920729252 & 8936528    &8936528.4  &0.6079270748 \\ 
14800000 &  5802667    &5802678.9 &0.3920720946 & 8997333    &8997321.1  &0.6079279054 \\ 
14900000 &  5841894    &5841886.2 &0.3920734228 & 9058106    &9058113.8  &0.6079265772 \\ 
15000000 &  5881111    &5881093.5 &0.3920740667 & 9118889    &9118906.5  &0.6079259333 \\ 
15100000 &  5920310    &5920300.8 &0.3920735099 & 9179690    &9179699.2  &0.6079264901 \\ 
15200000 &  5959515    &5959508.1 &0.3920733553 & 9240485    &9240491.9  &0.6079266447 \\ 
15300000 &  5998734    &5998715.3 &0.3920741176 & 9301266    &9301284.7  &0.6079258824 \\ 
15400000 &  6037936    &6037922.6 &0.3920737662 & 9362064    &9362077.4  &0.6079262338 \\ 
15500000 &  6077155    &6077129.9 &0.3920745161 & 9422845    &9422870.1  &0.6079254839 \\ 
15600000 &  6116369    &6116337.2 &0.3920749359 & 9483631    &9483662.8  &0.6079250641 \\ 
15700000 &  6155574    &6155544.5 &0.3920747771 & 9544426    &9544455.5  &0.6079252229 \\ 
15800000 &  6194781    &6194751.8 &0.3920747468 & 9605219    &9605248.2  &0.6079252532 \\ 
15900000 &  6233992    &6233959.1 &0.3920749686 & 9666008    &9666040.9  &0.6079250314 \\ 
16000000 &  6273204    &6273166.4 &0.3920752500 & 9726796    &9726833.6  &0.6079247500 \\ 
16100000 &  6312405    &6312373.7 &0.3920748447 & 9787595    &9787626.3  &0.6079251553 \\ 
16200000 &  6351612    &6351580.9 &0.3920748148 & 9848388    &9848419.1  &0.6079251852 \\ 
16300000 &  6390804    &6390788.2 &0.3920738650 & 9909196    &9909211.8  &0.6079261350 \\ 
16400000 &  6430021    &6429995.5 &0.3920744512 & 9969979    &9970004.5  &0.6079255488 \\ 
16500000 &  6469229    &6469202.8 &0.3920744848 & 10030771   &10030797.2 &0.6079255152 \\ 
16600000 &  6508430    &6508410.1 &0.3920740964 & 10091570   &10091589.9 &0.6079259036 \\ 
16700000 &  6547643    &6547617.4 &0.3920744311 & 10152357   &10152382.6 &0.6079255689 \\ 
16800000 &  6586845    &6586824.7 &0.3920741071 & 10213155   &10213175.3 &0.6079258929 \\ 
16900000 &  6626056    &6626032.0 &0.3920743195 & 10273944   &10273968.0 &0.6079256805 \\ 
17000000 &  6665261    &6665239.3 &0.3920741765 & 10334739   &10334760.7 &0.6079258235 \\ 
17100000 &  6704482    &6704446.6 &0.3920749708 & 10395518   &10395553.4 &0.6079250292 \\ 
17200000 &  6743707    &6743653.8 &0.3920759884 & 10456293   &10456346.2 &0.6079240116 \\ 
17300000 &  6782897    &6782861.1 &0.3920749711 & 10517103   &10517138.9 &0.6079250289 \\ 
17400000 &  6822092    &6822068.4 &0.3920742529 & 10577908   &10577931.6 &0.6079257471 \\ 
17500000 &  6861299    &6861275.7 &0.3920742286 & 10638701   &10638724.3 &0.6079257714 \\ 
17600000 &  6900508    &6900483.0 &0.3920743182 & 10699492   &10699517.0 &0.6079256818 \\ 
17700000 &  6939718    &6939690.3 &0.3920744633 & 10760282   &10760309.7 &0.6079255367 \\ 
17800000 &  6978917    &6978897.6 &0.3920739888 & 10821083   &10821102.4 &0.6079260112 \\ 
17900000 &  7018130    &7018104.9 &0.3920743017 & 10881870   &10881895.1 &0.6079256983 \\ 
18000000 &  7057349    &7057312.2 &0.3920749444 & 10942651   &10942687.8 &0.6079250556 \\ 
18100000 &  7096550    &7096519.5 &0.3920745856 & 11003450   &11003480.5 &0.6079254144 \\ 
18200000 &  7135765    &7135726.7 &0.3920750000 & 11064235   &11064273.3 &0.6079250000 \\ 
18300000 &  7174969    &7174934.0 &0.3920748087 & 11125031   &11125066.0 &0.6079251913 \\ 
18400000 &  7214172    &7214141.3 &0.3920745652 & 11185828   &11185858.7 &0.6079254348 \\ 
18500000 &  7253379    &7253348.6 &0.3920745405 & 11246621   &11246651.4 &0.6079254595 \\ 
18600000 &  7292578    &7292555.9 &0.3920740860 & 11307422   &11307444.1 &0.6079259140 \\ 
18700000 &  7331772    &7331763.2 &0.3920733690 & 11368228   &11368236.8 &0.6079266310 \\ 
18800000 &  7370984    &7370970.5 &0.3920736170 & 11429016   &11429029.5 &0.6079263830 \\ 
18900000 &  7410176    &7410177.8 &0.3920728042 & 11489824   &11489822.2 &0.6079271958 \\ 
19000000 &  7449388    &7449385.1 &0.3920730526 & 11550612   &11550614.9 &0.6079269474 \\ 
19100000 &  7488582    &7488592.4 &0.3920723560 & 11611418   &11611407.6 &0.6079276440 \\ 
19200000 &  7527801    &7527799.6 &0.3920729687 & 11672199   &11672200.4 &0.6079270313 \\ 
19300000 &  7566999    &7567006.9 &0.3920724870 & 11733001   &11732993.1 &0.6079275130 \\ 
19400000 &  7606189    &7606214.2 &0.3920715979 & 11793811   &11793785.8 &0.6079284021 \\ 
19500000 &  7645413    &7645421.5 &0.3920724615 & 11854587   &11854578.5 &0.6079275385 \\ 
19600000 &  7684616    &7684628.8 &0.3920722449 & 11915384   &11915371.2 &0.6079277551 \\ 
19700000 &  7723829    &7723836.1 &0.3920725381 & 11976171   &11976163.9 &0.6079274619 \\ 
19800000 &  7763022    &7763043.4 &0.3920718182 & 12036978   &12036956.6 &0.6079281818 \\ 
19900000 &  7802224    &7802250.7 &0.3920715578 & 12097776   &12097749.3 &0.6079284422 \\ 
20000000 &  7841425    &7841458.0 &0.3920712500 & 12158575   &12158542.0 &0.6079287500 \\ 
    \hline
 \caption*{\footnotesize Note:  $N$: Length of the block from $\mu(1)$ to $\mu(N)$; $N_{0}$: Frequency of $\mu(n)=0$ ;  $N_{0}^{\mbox{\tiny T}}(=N\times p_{_{\tiny{t}}})$: Frequency of $\mu(n)=0$ from number theory; $p_{_{\tiny{e}}}(=\frac{N_{0}}{N})$: Empirical probability; The number in bracket is theoretical probability from number theory $p_{_{\tiny{t}}}$. $N_{\vert\mu(n)\vert=1}$ and $N_{\vert\mu(n)\vert=1}^{\mbox{\tiny T}}$ are those for $\vert\mu(n)\vert=1$.} 
 \label{tb4}
\end{longtable}%


\begin{thebibliography}{99}
\bibitem{}
Chen NX, A new method for inverse black body radiation problem. Chinese Phys. Lett., 1987, 4: 337-340
\bibitem{}
Chen NX, A new method to introduce the M$\ddot{\mathrm{o}}$bius function and the M$\ddot{\mathrm{o}}$bius transform, J. Math. Phys., 1994, 35(6):3099-3108,doi: 10.1063/1.530455
\bibitem{} 
Chen NX, Chen Y, and Li GY, Theoretical investigation on inversion for the phonon density of states, Physics Letters A, 1990: 357-364
\bibitem{}
 Chen NX, M$\ddot{\mathrm{o}}$bius inversion in physics, Singapore: World Scientific, 2010
\bibitem{}
Chen NX, Modified M$\ddot{\mathrm{o}}$bius inverse formula and its applications in physics, Physical Review Letters, 1990, 64(11): 1193-1195
\bibitem{}
Good I J, Churchhouse R F, The Riemann hypothesis and pseudorandom features of the M$\ddot{\mathrm{o}}$bius sequence, Mathematics of Computation, 1968, 22(104): 857-861
\bibitem{}
Hardy G H, Wright E M. An introduction to the theory of numbers (Sixth Edition). Oxford University Press, 2008
\bibitem{}
Ji F M, Ye J P, Sun L, et al., An inverse transmissivity problem, its M$\ddot{\mathrm{o}}$bius inversion solution and new practical solution method. Phys Lett A, 2006, 352: 467-472
\bibitem{}
Kuznetsov E, Computing the Mertens function on a GPU,arXiv:1108.0135, 2011
\bibitem{}
Lioen WM, van de Lune J, Systematic Computations on Mertens' Conjecture and Dirichlet's Divisor Problem by Vectorized Sieving. In: From Universal Morphisms to Megabytes: A Baayen Space Odyssey. On the Occasion of the Retirement of P. C. Baayen (Ed. K. Apt, L. Schrijver, and N. Temme). Amsterdam, Netherlands: Stichting Mathematisch Centrum, Centrum voor Wiskunde en Informatica, 1994, pp. 421-432.
\bibitem{}
Odlyzko A M, te Riele H J J. Disproof of the Mertens conjecture, J. reine angew. Math, 1985, 357: 138-160
\bibitem{}
Sarnak P, Randomness in Number Theory, Asia Pacific Mathematics Newsletter, 2012, 2(3): 15-19
\bibitem{}
Schroeder M, Number theory in science and communication: with applications in cryptography, physics, digital information, computing, and self-similarity. Springer Science \& Business Media, 2008
\bibitem{}
Wang HY, Mathematical methods in Physics, Science Press, Beijing, 2013

\end{thebibliography}
\end{document}